\documentclass[lettersize,journal]{IEEEtran}
\usepackage{amsmath,amsfonts}
\usepackage{algorithmic}
\usepackage{algorithm}
\usepackage{array}
\usepackage[
caption=false,
font=normalsize,
labelfont=sf,
textfont=sf,
labelformat=simple,  % 使用简单格式
labelsep=space       % 编号与字母间加空格
]{subfig}

\usepackage{textcomp}
\usepackage{stfloats}
\usepackage{url}
\usepackage{verbatim}
\usepackage{graphicx}
\usepackage{amsmath}
\usepackage{mathrsfs}
\usepackage{epstopdf}
\usepackage{amssymb}
\usepackage{mathtools}

\usepackage[colorlinks,
linkcolor=blue,
anchorcolor=blue,
citecolor=blue,
urlcolor=blue%%设置DOI的颜色为蓝色
]{hyperref}

\newtheorem{remark}{Remark}
\newtheorem{assumption}{Assumption}

\newcommand{\ba}{\begin{array}}
	\newcommand{\ea}{\end{array}}
\newcommand{\be}{\begin{equation}}
	\newcommand{\ee}{\end{equation}}
\newcommand{\ben}{\begin{equation*}}
	\newcommand{\een}{\end{equation*}}
\newcommand{\bea}{\begin{eqnarray}}
	\newcommand{\eea}{\end{eqnarray}}
\newcommand{\bean}{\begin{eqnarray*}}
	\newcommand{\eean}{\end{eqnarray*}}
\newcommand{\beann}{\begin{align}}
	\newcommand{\eeann}{\end{align}}

\begin{document}
	
	\title{Periodic event-triggered impulsive control for fully heterogeneous stochastic multi-agent systems with a time-varying topology }
	
	\author{Xuetao Yang,~Ruilu An,~Quanxin Zhu$^{*}$,~\IEEEmembership{Senior Member,~IEEE}
	\thanks {	This work was jointly supported by the National Natural Science Foundation of China (62173139).}
		\thanks{Xuetao Yang is with School of Science, Nanjing University of Posts and Telecommunications, Nanjing 210023, Jiangsu, China, and also with MOE-LCSM, School of Mathematics and Statistics, Hunan Normal University, Changsha 410081, Hunan, China.}% <-this % stops a space
		\thanks{Ruilu An is with School of Science, Nanjing University of Posts and Telecommunications, Nanjing 210023, Jiangsu, China.}
		\thanks{Quanxin Zhu is with MOE-LCSM, School of Mathematics and
			Statistics, Hunan Normal University, Changsha 410081, China (The corresponding author is Quanxin Zhu with email:
			zqx22@126.com).}}
	
	% The paper headers
	%	\markboth{Journal of \LaTeX\ Class Files,~Vol.~14, No.~8, August~2021}%
	%	{Shell \MakeLowercase{\textit{et al.}}: A Sample Article Using IEEEtran.cls for IEEE Journals}

	\maketitle
	
	\begin{abstract}
		In this paper, we focus on a periodic event-triggered impulsive control (PETIC) for fully heterogeneous stochastic multi-agent systems (MASs) with a time-varying topology. Firstly, a novel time-varying topology is established by incorporating the energy consumption of each agent. This topology enables active adjustment of the information interaction intensity between agents. Secondly, to address the difficulties that agents with different dimensions cannot communicate in fully heterogeneous stochastic MASs, a virtual state space is designed. According to the above framework, novel PETICs with/without actuation delays are presented to achieve the mean-square exponential consensus of fully heterogeneous stochastic MASs. Finally, the effectiveness of the proposed methods is verified through a numerical simulation of unmanned aerial vehicles and unmanned ground vehicles.
	\end{abstract}
	
	
	\section{Introduction}
	
	\IEEEPARstart{I}{n} recent years, multi-agent systems (MASs) have attracted remarkable attention due to their extensive applications in engineering fields such as unmanned aerial vehicle motion, smart grids, robot control and so on (see, e.g., \cite{LTZY24, ZYCSXD25, WWZVC21}). One of the core issues underlying the effective operation of multi-agent systems is the consensus problem, which enables all agents to achieve a common task objective by receiving information from their adjacent agents (see, e.g., \cite{ZLZ25, ZCDX24}). Due to the importance of the consensus problem, it has gradually become a popular topic, which promotes our research.
	
	In engineering practice, the condition that all agents have identical models exhibits significant limitations, which has driven the development of heterogeneous MASs. It is worth noting that most heterogeneous MASs primarily differ in the rules governing each agent's state evolution, while differences in  agents' dimensions are often overlooked (see, e.g., \cite{SLWGW23, YJLZY24, WD24, GLZA23, WDLC24}). This observation has aroused our research interest in fully heterogeneous MASs with agents of different dimensions. However, differences in agents' dimensions create obstacles to information transmission between agents. Therefore, it is a challenge to achieve consensus in fully heterogeneous MASs. In \cite{HZWR23}, a dynamic compensator is constructed for each agent to receive and transmit information. In \cite{LXLF22}, a distributed observer is proposed to estimate the states of follower systems and to ensure information transmission between followers. However, the aforementioned fully heterogeneous MASs are all deterministic systems. For fully heterogeneous stochastic MASs, the stochastic noise is unobservable, which leads to the methods applicable to deterministic systems infeasible. Thus, for fully heterogeneous stochastic MASs, how to enable agents to overcome dimensional differences and achieve consensus is one of the key problems addressed in this paper.
	
	Currently, the topological structures in most studies are static, which means that fixed communication is always maintained between agents (see, e.g., \cite{KRSSYD24, YLCLY25, HZH25, LSL24}). However, in reality, due to factors such as environmental changes, the agents' own conditions and task requirements, the information interaction between agents changes over time. In dynamic topological scenarios, controllers designed for static topologies often lose their effectiveness because they cannot adapt to the time-varying communication relationships between agents. To address this issue, some studies have proposed controllers based on switching topologies (see, e.g., \cite{ZZXX24, WWD22, FHWWY23}). Nevertheless, the switching modes employed in these methods are limited and preconfigured. In this paper, we are devoted to designing a time-varying topology that adapts to the individual characteristics of agents, which facilitates MASs in achieving its task objectives. Energy is a fundamental requirement for the stable operation of agents (see, e.g., \cite{BKAT22, YCY20, YBZ22}). When energy is insufficient, the risk of agent instability increases significantly. Therefore, it is essential to design an energy-adaptive dynamic topological structure. Then, how to design a control strategy with an energy-adaptive dynamic topology has become an urgent problem that  needs solving.
	
	During the past several years, impulsive controls have attracted widespread attention due to its characteristics of simple design and easy implementation. Traditional impulsive controls are triggered based on time (see, e.g., \cite{HM23, LXS25, LJZ22}). However, event-triggered impulsive controls (ETICs) can provide flexible updating times for the controller based on an event-triggering mechanism (ETM) (see, e.g., \cite{JWLG24, PZLZ25, ZPLWC22, ZB23, HM22}). So, ETICs can significantly reduce bandwidth consumption and they have been extensively applied in many systems such as complex networks, switched systems and coupled neural networks (see, e.g., \cite{ZCWH25, RFXWZ22, ZCL24, SYLHC22}). Therefore, we are committed to applying ETICs to fully heterogeneous stochastic MASs. 

In the scenario of time-varying topology, we combine an energy consumption function
with periodic ETM (PETM) and propose a novel periodic ETIC (PETIC). The advantage is that this novel PETIC can not only solve the
consensus problem of
 stochastic MASs, but also  avoid the Zeno phenomenon effectively. Meanwhile, to account for the impact of actuation delays (i.e., the delays occurring between the controller and the actuator), we also propose a novel PETIC with actuation delays. 
The significant contributions of this paper are summarized as follows:
\begin{enumerate}
	\item{To address the issue of different states' dimensions in fully heterogeneous stochastic MASs, a novel virtual state space is introduced, which ensures the transmission of information between agents. Moreover, the virtual system states can be effectively degenerated into the original system states.}
	
	\item{Different from the topologies proposed in \cite{KRSSYD24, YLCLY25, HZH25, LSL24, ZZXX24, WWD22, FHWWY23}, we establish a novel time-varying topology based on energy. By incorporating the energy consumption of each agent into the topology, the information transmission between adjacent agents is more adaptive to the individual characteristics of each agent.}
	
	\item {Two kinds of novel PETICs are proposed. Different from the fixed updating instants of impulsive controls (see, e.g., \cite{HM23, LXS25, LJZ22}), control updating instants considered in this paper are a sequence of stopping times selected by PETM, which save communication resources more reasonably. Meanwhile, unlike traditional PETMs 
based on system states directly (see, e.g., \cite{JWLG24, PZLZ25, HM22}), our PETM is based on virtual states, which is more suitable for fully heterogeneous systems. Furthermore, energy consumption designed in PETICs plays a positive role in accelerating the consensus rate of stochastic MASs.}

\end{enumerate}

This paper is structured as follows.  In Section II, we introduce a fully heterogeneous stochastic MAS with a time-varying topology.  In Section III, we propose PETICs with/without actuation delay strategies to ensure system consensus.  In Section IV, we provide the simulation results for unmanned aerial vehicles (UAVs) and unmanned ground vehicles (UGVs).  In the final section, we wrap up the paper with some conclusions.

\section{Problem formulation }

\subsection{Notations}
Let $\mathbb{N}=\lbrace0,1,\cdots\rbrace$ and  $\mathbb{N}^+=\lbrace1,2,\cdots\rbrace$. For a real square matrix $M$, $\lambda_{\max}(M)$ and $\lambda_{\min}(M)$  denote the largest and smallest eigenvalues, respectively. $|M|$ represents the Euclidean norm. For real square matrices $M_1$ and $M_2$, $M_1\otimes M_2$ is the Kronecker product of the matrices $M_1$ and $M_2$. For constants $a,b$ and $c$, $diag\lbrace a,b,c\rbrace$ denotes a square matrix with diagonal elements $a, b$ and $c$, while all other elements are zero. Furthermore, define $a\vee b=\max\{a,b\}$. $I_n$ refers to the identity matrix of dimension $n$. $\mathbf{0}$ is a dimension-compatible
zero matrix. $\omega(t)$ is a standard Wiener process defined on a complete probability space $(\Omega,\mathcal{F},\lbrace\mathcal{F}_t\rbrace_{t\geq 0},\mathcal{P})$, in which the filter $\lbrace\mathcal{F}_t\rbrace_{t\geq 0}$ satisfies the usual conditions.

\subsection{Time-varying Graph Theory }

For a positive integer $N$, let $\mathcal{N}=\{1,2,\cdots,N\}$. In this paper, $\tau_i(t)$, $i\in \mathcal{N}$ represents the energy consumption of each agent (see \cite{BKAT22} for instance). We incorporate $\tau_i(t)$ into the communication topology. Then, the communication topology of followers is described as a directed graph $\mathcal{G}(t)=(\mathbb{V},\mathcal{E}(t),\mathcal{A}(t))$ with the set $\mathbb{V}=\{v_1,v_2,\cdots,v_{N}\}$, the edge set $\mathcal{E}(t)\subseteq \mathbb{V}\times\mathbb{V}$ and the weighted adjacency matrix $\mathcal{A}(t)=[a_{ij}(t)]_{N\times N}$, where $a_{ij}(t)=\bar{a}_{ij}e^{-\alpha\tau_{i}(t)}$, $\bar{a}_{ij}$ and $\alpha$ are positive constants. The neighbors of agent
$i$ at instant $t$ are depicted as $\mathcal{N}_i(t)=\{j\in\mathbb{V}|(j,i)\in\mathcal{E}(t)\}$. The Laplacian-like matrix of $\mathcal{G}(t)$ is given as $\mathbb{L}=[\bar{a}_{ij}]_{N\times N}$ with $l_{ii}=\Sigma_{i\neq j}\bar{a}_{ij}$ and $l_{ij}=-\bar{a}_{ij},~i\neq j$.
If the matrix $\mathcal{G}(t)$ has a leader, then the diagonal matrix $\mathcal{B}(t)=[b_i(t)]_{N\times N}=diag\{\bar{b}_1e^{-\alpha\tau_1(t)},~\bar{b}_2e^{-\alpha\tau_2(t)},\cdots,\bar{b}_{N}e^{-\alpha\tau_{N}(t)}\}$ can be used to represent the connectivity between leaders and followers, where $\bar{b}_i$ is a non-negative constant. We define a matrix $B=diag\{\bar{b}_1,\bar{b}_2,\cdots,\bar{b}_{N}\}$ and an information-exchange matrix $H=-\mathcal{L}+B$, where $B$ is not a zero matrix.
\begin{assumption}
	Suppose that the energy consumption $\tau_i(t)$ is independent of the system states $x_0(t)$, $x_i(t)$, $z_i(t)$ and $y_i(t)$, which will be determined later by (\ref{01}), (\ref{02}), (\ref{03}), and (\ref{05}).
\end{assumption}

\begin{remark}
	Different from the traditional static communication topologies (see, e.g., \cite{KRSSYD24, YLCLY25, HZH25, LSL24}), this paper explicitly considers the influence of energy consumption on information interaction between adjacent agents. To capture this impact, the novel time-varying weighted adjacency matrix $\mathcal{A}(t)$ and a diagonal matrix $\mathcal{B}(t)$ are presented.
	Then, topology structure under different energy consumption scenarios is as follows:
	\begin{enumerate}
		\item{When $\tau_i(t) = 0$, the weights $a_{ij}(t) = \bar{a}_{ij}$ and $b_i(t)=\bar{b}_i$. That is to say, when an agent has sufficient energy, the topology assigns the agent the maximum communication weight. At this point, $\mathcal{G}(t)$ degenerates into a static topology $\mathcal{G}$, which is consistent with traditional communication graphs.}
		\item{When $0 < \tau_i(t) < +\infty$, the weights $a_{ij}(t)$ and $b_i(t)$ decrease as $\tau_i(t)$ increases. The more energy an agent consumes, the less capability it has to receive information transmitted by other agents. This dynamic variation also prevents high-energy-consumption agents from undertaking excessive communication and coordination tasks.}
		\item{When $\tau_i(t)\rightarrow +\infty$, the weights $a_{ij}(t) \rightarrow 0$ and $b_i(t)\rightarrow 0$. That is, when an agent runs out of energy, it is automatically removed from the topological structure. }
	\end{enumerate}
\end{remark}
\subsection{A fully heterogeneous stochastic MAS}
Take into account the leader state below
\begin{eqnarray}
	&&dx_0(t)=[C_0x_0(t)+f_{x0}(x_0(t))]dt\nonumber\\
	&&\quad\quad\quad\quad+D_0x_0(t)d\omega(t),~t> 0,\label{01}
\end{eqnarray}
where $x_0(t)\in \mathbb{R}^{n_0}$ is the leader state with an initial value $x_0(0)\in \mathbb{R}^{n_0}$. $C_0\in \mathbb{R}^{n_0\times n_0}$ and $D_0\in  \mathbb{R}^{n_0\times n_0}$ are constant matrices. $f_{x0}(\cdot): \mathbb{R}^{n_0}\rightarrow \mathbb{R}^{n_0}$ is a continuous function.

Meanwhile, consider the following stochastic MAS of $N$ agents:
\begin{eqnarray}
	&&dx_i(t)=[C_ix_i(t)+f_{xi}(x_i(t))+u_i(t)]dt\nonumber\\
	&&\quad\quad\quad\quad+D_ix_i(t)d\omega(t),~i\in \mathcal{N},~t> 0,\label{02}
\end{eqnarray}
where $x_i(t)\in \mathbb{R}^{n_i}$ represents the i-th system state with an initial value $x_i(0)\in \mathbb{R}^{n_i}$. $u_i(t)\in \mathbb{R}^{n_i}$ stands for a control input. $C_i\in \mathbb{R}^{{n_i}\times {n_i}}$ and $D_i\in \mathbb{R}^{{n_i}\times{n_i}}$ are system matrices. $f_{xi}(\cdot): \mathbb{R}^{n_i}\rightarrow \mathbb{R}^{n_i}$ is a continuous function. Next, to solve the consensus problem of the leader-following system, we propose the following assumption.

\begin{assumption}
	Suppose that there exists a matrix $\Xi_i\in \mathbb{R}^{n_i\times n_0}$ such that $\Xi_iC_0=C_i\Xi_i$ and $\Xi_iD_0=D_i\Xi_i$.
\end{assumption}
\begin{remark}
	Assumption 1 is a classical condition for heterogeneous leader-following systems (see, e.g., \cite{HZWR23, LXLF22}).
\end{remark}

Let $z_i(t)=x_i(t)-\Xi_ix_0(t)$ with the initial value $z_i(0)=x_i(0)-\Xi_ix_0(0)$. Under Assumption 2, systems (\ref{01}) and (\ref{02}) can be rewritten as
\begin{eqnarray}
	&&dz_i(t)=[C_iz_i(t)+f_i(z_i(t))+u_i(t)]dt\nonumber\\
	&&\quad\quad\quad\quad+D_iz_i(t)d\omega(t),~i\in \mathcal{N},~t> 0,\label{03}
\end{eqnarray}
where $f_i(z_i(t))=f_{xi}(x_i(t))-\Xi_if_{x0}(x_0(t))$. Subsequently, one assumption regarding system (\ref{03}) are presented.

\begin{assumption}
	Suppose that $f_{x0}(x_0(t))$, $f_{xi}(x_i(t))$ and $f_i(z_i(t))$ satisfy the Lipschitz condition with $f_{x0}(0)=0$, $f_{xi}(0)=0$ and $f_i(0)=0$. Moreover, there exists a nonnegative constant $L_{f_i}$ such that $|f_i(z_i(t))|\leq L_{f_i}|z_i(t)|$.
\end{assumption}

\begin{remark}
	Under Assumption 3, it follows from \cite{YZW24} that system (\ref{03}) has a unique global solution.
\end{remark}

\subsection{Virtual state space}
To solve the problem of different state dimensions in the fully heterogeneous stochastic MAS, we construct a state space with $m$ dimensions. Then, the virtual state $y_i(t)$ satisfies
\begin{eqnarray}
	y_i(t)=\Phi_iz_i(t),~i\in\mathcal{N},\label{04}
\end{eqnarray}
where $\Phi_i\in \mathbb{R}^{m\times n_i}$ is a constant matrix and $m\geq n_i$.   In the following, we present a necessary assumption to guarantee the effectiveness of the virtual state space.
\begin{assumption}
	Suppose that matrix $\Phi_i$ satisfies $rank(\Phi_i)=n_i$. Moreover, there exists a matrix $\Theta_i\in \mathbb{R}^{n_i\times m}$ such that $\Theta_i\Phi_i=I_{n_i}$.
\end{assumption}
\begin{remark}
	Under Assumption 4, it follows from formula (\ref{04}) that $z_i(t)=\Theta_iy_i(t)$. The existence of matrix $\Theta_i$ implies that the original state $z_i(t)$ can be exactly reconstructed from the virtual state $y_i(t)$ via a fixed linear transformation. Consequently, embedding the original states of heterogeneous systems into an m-dimensional virtual state space does not lead to any loss of state information.
\end{remark}
\begin{remark}
	The necessity for an m-dimensional virtual state space to be constructed is stated as follows: On the one hand, compared with \cite{HZWR23} and \cite{LXLF22}, system (\ref{03}) includes unobservable stochastic noise. On the other hand, the condition $m\geq n_i$ is required to preserve all information in the following systems. If each follower is directly mapped to the leader's dimension, it cannot cover the application scenarios when $n_0<n_i$. Therefore, we design the virtual state space, which not only eliminates the need to observe stochastic noise but also covers all application scenarios of leader-follower systems.
\end{remark}

According to system (\ref{03}) and formula (\ref{04}), we obtain the virtual stochastic MAS
\begin{eqnarray}
	&&dy_i(t)=[\Phi_iC_i\Theta_iy_i(t)+\Phi_if_i(t)+\Phi_iu_i(t)]dt\nonumber\\
	&&\quad\quad\quad\quad+\Phi_iD_i\Theta_iy_i(t)d\omega(t),~i\in \mathcal{N},~t>0,\label{05}
\end{eqnarray}
where $f_i(t)=f_i(\Theta_iy_i(t))$ and $y_i(0)=\Phi_iz_i(0)$ denotes the initial value of the state $y_i(t)$.

\section{Main results}\label{sec-3}
\begin{figure}[H]
\centering
\includegraphics[width=9cm]{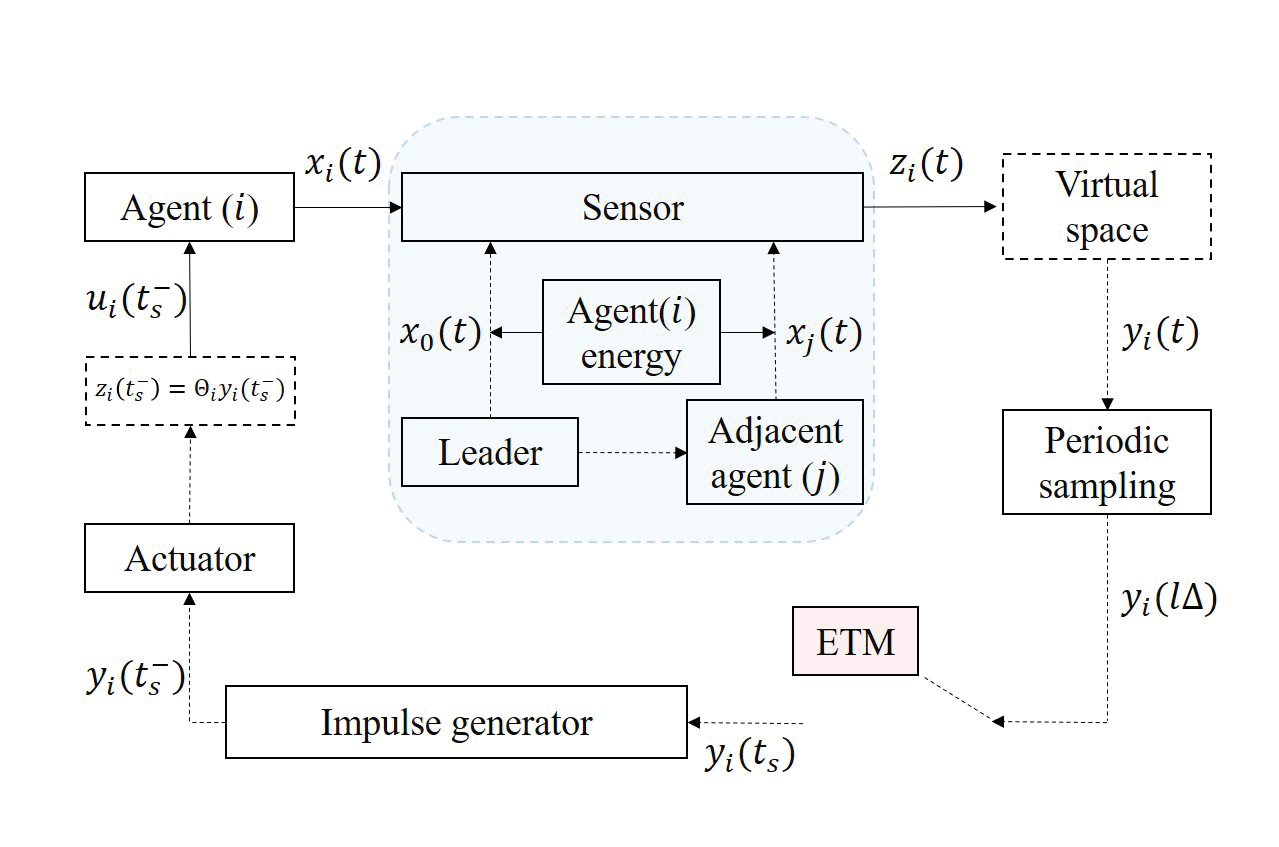}
\caption{ Configuration of a stochastic MAS with time-varying topology.}
\label{fig0}
\end{figure}

In this section, we investigate PETIC for a stochastic MAS with time-varying topology, as illustrated in Fig. \ref{fig0}. Specifically, we consider two cases: a controller without actuation delays and a controller with actuation delays.
\vspace{-1em}

\subsection{PETIC}
In order to reduce the consumption of control resources, the impulsive controller is designed as
\begin{eqnarray}
	&&u_i(t)=\Sigma_{s=0}^\infty k_iU_i(t)\delta(t-t_s),\label{00}
\end{eqnarray}
where  $k_i$ is a negative constant and $U_i(t)=e^{-\alpha\tau_i(t)}[\Sigma_{j= 1}^N(\bar{a}_{ij}\\\times\Theta_i(y_i(t)-y_j(t)))+\bar{b}_i\Theta_iy_i(t)]$. As defined in \cite{HM23}, $\delta(\cdot)$ represents the Dirac function. ${\lbrace t_s\rbrace}_{s\in\mathbb{N}}$ is an increasing (stopping time) impulsive sequence with $t_0\neq0$ and it will be determined later by (\ref{007}) and (\ref{07}). $y_i(t)$ is right continuous with left limits at impulsive instants $t_s$, i.e., $y_i(t_s)=y_i(t_s^+)$.
\begin{remark}
	Compared to controllers that rely on the system's original states (see, e.g., \cite{PZLZ25, ZPLWC22}), PETIC in this paper is based on the virtual state $y_i(t)$ with the same dimension, which facilitates the comparison of states between adjacent agents. Additionally, the energy consumption $\tau_i(t)$ is incorporated into PETIC. Due to the exponential decay property of $e^{-\alpha\tau_i(t)}$, the system can achieve stability more easily.
\end{remark}

Let $y(t)=[y_1^\mathrm{T}(t),y_2^\mathrm{T}(t),\cdots,y_N^\mathrm{T}(t)]^\mathrm{T}$. Then, stochastic MAS (\ref{05}) is derived as
\begin{eqnarray}
	\left\{\begin{aligned}
		&dy(t)=[\Phi C\Theta y(t)+\Phi f(t)]dt\\
		&\quad\quad\quad~+\Phi D\Theta y(t)d\omega(t),~t\neq t_s,\\
		&y(t)=e^{-\alpha\Gamma(t)}K\tilde{H}\Phi\Theta y(t^-)+y(t^-),~t = t_s,
	\end{aligned}\label{06}
	\right.
\end{eqnarray}
where
\begin{eqnarray*}
	&&f(t)=[f_1^\mathrm{T}(t),f_2^\mathrm{T}(t),\cdots,f_N^\mathrm{T}(t)]^\mathrm{T},\\
	&&\Phi=diag\{\Phi_1,\Phi_2,\cdots,\Phi_N\},\\
	&&\Theta=diag\{\Theta_1,\Theta_2,\cdots,\Theta_N\},\\
	&&C=diag\{C_1,C_2,\cdots,C_N\},\\
	&&D=diag\{D_1,D_2,\cdots,D_N\},\\
	&&\tilde{H}=I_m\otimes H,~K=I_m\otimes diag\{k_1,k_2,\cdots,k_N\},\\
	&&\Gamma(t)=I_m\otimes diag\{\tau_1(t),\tau_2(t),\cdots\tau_N(t)\}.
\end{eqnarray*}

Given the sampling period $\Delta$, PETM is designed as
\begin{align}
		&t_{0}=\mathop{\min}_{\iota\in \mathbb{N}^+}\{\iota\Delta:W(y(\iota\Delta))>\psi_1 W(y(0))\},\label{007}\\
		&t_{s+1}=t_s+\mathop{\min}_{\iota\in \mathbb{N}^+}\{\iota\Delta:W(y(t_s+ \iota\Delta)>\psi_2 W(y(t_s))\},
		\label{07}
\end{align}
where $W(y(t))=e^{\gamma t}y^\mathrm{T}(t)Py(t)$, constants $\gamma>0$, $\psi_1\geq1$, $\psi_2\geq1$ and $P\in \mathbb{R}^{m N\times mN}$ is a positive definite matrix.
\begin{remark}
	Compared to the traditional PETMs in \cite{JWLG24, PZLZ25, HM22}, this paper takes into account systems' virtual states $y(t)$. Additionally, we adopt the periodic observations to generate a minimum inter-event time $\Delta$, which avoids the Zeno phenomenon.
\end{remark}

In what follows, we propose two necessary assumptions that play a crucial role in the proof of Theorem 1.

\begin{assumption}
	Suppose that the energy consumption $\tau_{i}(t)$ satisfies $\tau_{i}(t)\geq \tau$.
\end{assumption}
\begin{remark}
In Assumption 5, $\tau$ represents the minimum energy consumption of the agent. Under PETIC (\ref{00}), the agent can achieve stability only when its energy consumption is higher than the minimum energy consumption. Under Assumption 5, we have $e^{-\alpha\tau_i(t)}\leq e^{-\alpha\tau}$, which is necessary for the proof of Theorem 1.
\end{remark}

\begin{assumption}\label{A3}
	Suppose that
	\begin{eqnarray*}
		2\alpha\tau-ln\psi_2-ln\lambda_1-\lambda\Delta>0,
	\end{eqnarray*}
	where  $L_f=diag\{L_{f_1}I_m,L_{f_2} I_m,\cdots,L_{f_N}I_m\}$,  $\lambda=\max\{0,~\lambda_{\max}([P\Phi C\Theta+(\Phi C\Theta)^\mathrm{T}P+(\Phi D\Theta)^\mathrm{T}P\Phi D\Theta+P^\mathrm{T}P+(L_f\Theta)^\mathrm{T} L_f\Theta]P^{-1})\}$ and $\lambda_1=\lambda_{\max}((e^{\alpha\tau}I_{mN}+K\tilde{H}\Phi\Theta)^\mathrm{T}P(e^{\alpha\tau}I_{mN}+K\tilde{H}\Phi\Theta)P^{-1})$.
\end{assumption}
\begin{remark}
	Let $\bar{\gamma}=\frac{2\alpha\tau-ln\psi_2 -ln\lambda_1 -\lambda\Delta}{\Delta}$. Under Assumption 6, we know $\bar{\gamma}>0$. Choose a constant $\gamma\in(0,\bar{\gamma}]$, we have $
	\psi_2\lambda_1 e^{(\lambda+\gamma)\Delta-2\alpha\tau}\leq 1$, which is important for the proof of Theorem 1. Moreover, it is evident that if other parameters are fixed, then the larger the value of $\alpha$, the larger the value of $\bar{\gamma}$. In other words, $\alpha$ is useful for adjusting the convergence speed of the system.
\end{remark}

Now, we are ready to state the main results on the stabilization of system (\ref{06}) under PETIC.

Theorem 1: Under Assumptions 1-6, system (\ref{06}) satisfies $\mathbb{E}|y(t)|^2\leq Me^{-\gamma t}\mathbb{E}|y(0)|^2$, where $M$ and $\gamma$ are positive constants.

\noindent{\bf Proof.}
Select a Lyapunov function $V(y(t))=y^\mathrm{T}(t) Py(t)$. By It\^{o}'s formula and Assumption 3,  we have
\begin{eqnarray}
	&&\mathcal{L}V(y(t))\nonumber\\
	&&=y^\mathrm{T}(t)[P\Phi C\Theta+(\Phi C\Theta)^\mathrm{T}P\nonumber\\
	&&\quad+(\Phi D\Theta)^\mathrm{T}P\Phi D\Theta]y(t)+2y^\mathrm{T}(t)P\Phi f(t)\nonumber\\
	&&\leq y^\mathrm{T}(t)[P\Phi C\Theta+(\Phi C\Theta)^\mathrm{T}P\nonumber\\
	&&\quad+(\Phi D\Theta)^\mathrm{T}P\Phi D\Theta+P^\mathrm{T}P+(L_f\Theta)^\mathrm{T}L_f\Theta]y(t)\nonumber\\
	&&\leq\lambda V(y(t)).\label{10}
\end{eqnarray}

Since $W(y(t))=e^{\gamma t}y^\mathrm{T}(t)Py(t)$, we obtain
\begin{eqnarray*}
	\frac{d\mathbb{E}W(y(t))}{dt}= \mathbb{E}\mathcal{L}W(y(t))\leq(\lambda+\gamma)\mathbb{E}W(y(t)).
\end{eqnarray*}

For $t\in[t_s+l\Delta,t_s+(l+1)\Delta)$, $l\in\mathbb{N}$, using the comparison principle, we get
\begin{eqnarray}
	&&\mathbb{E}W(y(t))\nonumber\\
	&&\leq e^{(\lambda+\gamma)(t-t_s-l\Delta)}\mathbb{E}W(y(t_s+l\Delta))\nonumber\\
	&&\leq e^{(\lambda+\gamma)\Delta}\mathbb{E}W(y(t_s+l\Delta)).\label{11}
\end{eqnarray}

When $l=0$, applying system (\ref{06}) and Assumption 5, (\ref{11}) can be derived as
\begin{eqnarray}
	&&\mathbb{E}W(y(t))\nonumber\\
	&&\leq e^{(\lambda+\gamma)\Delta}\mathbb{E}W(y(t_s))\nonumber\\
	&&=e^{(\lambda+\gamma)\Delta+\gamma t_s}\mathbb{E}[y^\mathrm{T}(t_s^-)(I_{mN}+e^{-\alpha\Gamma(t_s)}K\tilde{H}\Phi\Theta)^\mathrm{T}P\nonumber\\
	&&\quad\times(I_{mN}+e^{-\alpha\Gamma(t_s)}K\tilde{H}\Phi\Theta)y(t_s^-)]\nonumber\\
	&&\leq \lambda_1 e^{(\lambda+\gamma)\Delta-2\alpha\tau}\mathbb{E}W(y(t_s^-)).\label{12}
\end{eqnarray}

When $l\in(0,l_s)$, $l_s=\max\{\iota\in \mathbb{N}^+: t_s+\iota\Delta<t_{s+1}\}$. According to PETM (\ref{07}) and Assumption 5, (\ref{11}) implies
\begin{eqnarray}
	&&\mathbb{E}W(y(t))\nonumber\\
	&&\leq e^{(\lambda+\gamma)\Delta}\mathbb{E}W(y(t_s+l\Delta))\nonumber\\
	&&\leq\psi_2e^{(\lambda+\gamma)\Delta}\mathbb{E}W(y(t_s))\nonumber\\
	&&\leq\psi_2\lambda_1 e^{(\lambda+\gamma)\Delta-2\alpha\tau}\mathbb{E}W(y(t_s^-)).\label{13}
\end{eqnarray}

When $t\in[t_s+l_s\Delta,t_{s+1})$, similar to (\ref{13}), we have
\begin{eqnarray}
	&&\mathbb{E}W(y(t))\nonumber\\
	&&\leq e^{(\lambda+\gamma)\Delta}\mathbb{E}W(y(t_s+l_s\Delta))\nonumber\\
	&&\leq\psi_2\lambda_1 e^{(\lambda+\gamma)\Delta-2\alpha\tau}\mathbb{E}W(y(t_s^-)).\label{013}
\end{eqnarray}

Since $\psi_2\geq 1$, it follows from (\ref{12}), (\ref{13}) and (\ref{013}) that for $t\in[t_s,t_{s+1})$,
\begin{eqnarray}
	&&\mathbb{E}W(y(t))\leq \psi_2\lambda_1 e^{(\lambda+\gamma)\Delta-2\alpha\tau}\mathbb{E}W(y(t_s^-)).\label{0013}
\end{eqnarray}

According to Assumption 6 and $\gamma\in(0,\bar{\gamma}]$, (\ref{0013}) implies $\mathbb{E}W(y(t))\leq\mathbb{E}W(y(t_s^-))$. Thus, for $t\geq t_0$, we obtain
\begin{eqnarray}
	&&\mathbb{E}W(y(t))\leq\mathbb{E}W(y(t_s^-))\leq\mathbb{E}W(y(t_{s-1}^-))\nonumber\\
	&&\quad\qquad\quad~\leq\cdots\leq\mathbb{E}W(y(t_0^-)).\label{14}
\end{eqnarray}

Then, for $t\in[t_0-\Delta,t_0)$, applying the comparison principle and PETM (\ref{007}), we can
reveal that
\begin{eqnarray}
	&&\mathbb{E}W(y(t))\nonumber\\
	&&\leq e^{(\lambda+\gamma)\Delta}\mathbb{E}W(y(t_0^--\Delta))\nonumber\\
	&&\leq\psi_1e^{(\lambda+\gamma)\Delta}\mathbb{E}W(y(0)).\label{0014}
\end{eqnarray}

Combining (\ref{14}) and (\ref{0014}), for $t\geq0$, we have $\mathbb{E}W(y(t))\leq\psi_1e^{(\lambda+\gamma)\Delta}\mathbb{E}W(y(0))$. It is evident that $\lambda_{\min}(P)|y(t)|^2\leq V(y(t))\leq\lambda_{\max}(P)|y(t)|^2$. Then, we know
\begin{eqnarray*}
	\mathbb{E}|y(t)|^2\leq Me^{-\gamma t}\mathbb{E}|y(0)|^2,
\end{eqnarray*}
where $M=\psi_1e^{(\lambda+\gamma)\Delta}\frac{\lambda_{\max}(P)}{\lambda_{\min}(P)}$. Therefore, system (\ref{06}) is mean-square exponentially stable. This completes the proof of Theorem 1.
{\hfill $\Box$}

\subsection{PETIC with actuation delay}

The controller with actuation delay is designed as
\begin{eqnarray}
	&&\tilde{u}_i(t)=\Sigma_{s=1}^\infty \tilde{k}_i\tilde{U}_i(t)\delta(t-t_s),\label{000}
\end{eqnarray}
where $\tilde{U}_i(t)=e^{-\alpha\tau_i(t)}[\Sigma_{j= 1}^N[\bar{a}_{ij}\Theta_i(y_i(t-\tau_s)-y_j(t-\tau_s))]+\bar{b}_i\Theta_iy_i(t-\tau_s)]$, $\tilde{k}_i\in \mathbb{R}$ is a constant and $\tau_s\ge0$ is the actuation delay.  ${\lbrace t_s\rbrace}_{s\in\mathbb{N}}$ is also determined by (\ref{007}) and (\ref{07}) with $t_0\neq 0$. Then, system (\ref{06}) can be rewritten as
\begin{eqnarray}
	\left\{\begin{aligned}
		&dy(t)=[\Phi C\Theta y(t)+\Phi f(t)]dt\\
		&\quad\quad\quad~+\Phi D\Theta y(t)d\omega(t),~t\neq t_s,\\
		&y(t)=e^{-\alpha\Gamma(t)}\tilde{K}\tilde{H}\Phi\Theta y(t^--\tau_s),~t = t_s,
	\end{aligned}\label{15}
	\right.
\end{eqnarray}
where $\tilde{K}=I_m\otimes diag\{\tilde{k}_1,\tilde{k}_2,\cdots,\tilde{k}_N\}$.

\begin{remark}
	Compared to PETIC (\ref{00}), the control input $\tilde{u}_i(t)$ of PETIC (\ref{000}) is constructed based on the virtual state $y_i(t-\tau_s)$ at the actual impulsive control execution instants. This approach prevents control ineffectiveness caused by the mismatch between the designed control gain and the actual system state. When the actuation delay $\tau_s = 0$, PETIC (\ref{000}) degenerates into PETIC (\ref{00}). Furthermore, unlike the control gain $k_i$ in PETIC (\ref{00}), which must be negative, the $\tilde{k}_i$ in PETIC (\ref{000}) is not subject to this strict constraint. This flexibility in gain selection expands the applicable scenarios of the control strategy.
\end{remark}

Next, we introduce three necessary assumptions that play a vital role in achieving the stabilization of system (\ref{15}).

\begin{assumption}
	Suppose that the energy consumption $\tau_{i}(t)$ satisfies $\tau_{i}(t)\geq \beta_i t$, where $\beta_i$ is a positive constant.
\end{assumption}
\begin{remark}
	Based on \cite{BKAT22}, $\beta_i$ denotes the minimum power of each agent. Different from PETIC (\ref{00}), which requires the energy consumption $\tau_i(t) \geq \tau$, under PETIC (\ref{000}), the system can achieve stability when the energy consumption satisfies $\tau_i(t) \geq 0$. Under Assumption 7, we have $e^{-\alpha\tau_i(t)}\leq e^{-\alpha\beta_i t}$, which is necessary for the proofs of Theorem 2.
\end{remark}
\begin{assumption}
	Suppose that the actuation delay $\tau_s$ satisfies $\tau_s<\Delta$, $s\in\mathbb{N}$.
\end{assumption}
\begin{remark}
	 Assumption 8 describes a typical scenario of control delays under PETC (see \cite{ZCWH25}). Under Assumption 8, actual control instants always lie between two adjacent triggering instants.
\end{remark}
\begin{assumption}\label{A6}
	Suppose that
	\begin{eqnarray*}
		(2\alpha\beta-\lambda)\Delta-ln\psi_2-ln\tilde{\lambda}_1>0,
	\end{eqnarray*}
\end{assumption}
where $\beta=\max\{\beta_1,\beta_2,\cdots,\beta_N\}$ and $\tilde{\lambda}_1=\lambda_{\max}((\tilde{K}\tilde{H}\Phi\Theta)^\mathrm{T}\\\times P(\tilde{K}\tilde{H}\Phi\Theta)P^{-1})$.
\begin{remark}
	Let $\bar{\gamma}^{'}=\frac{(2\alpha\beta-\lambda)\Delta-ln\psi_2-ln\tilde{\lambda}_1}{2\Delta}$. Based on Assumption 9, we know $\bar{\gamma}^{'}>0$. Choose a constant $\gamma\in(0,\bar{\gamma}^{'}]$, we have $
\psi_2\tilde{\lambda}_1e^{(\lambda+2\gamma-2\alpha\beta) \Delta}\leq 1$.
\end{remark}

Now, we are prepared to introduce the main results on the stabilization of system (\ref{15}) under PETIC with actuation delays.

Theorem 2: Under Assumptions 1-4 and 7-9, system (\ref{15}) satisfies $\mathbb{E}|y(t)|^2\leq Me^{-\gamma t}\mathbb{E}|y(0)|^2$, where $M$ and $\gamma$ are positive constants.

\noindent{\bf Proof.}
Similar to the proof of Theorem 1 for $t\in[t_s+l\Delta,t_s+(l+1)\Delta)$, it is easy to know $\mathbb{E}W(y(t))\leq e^{(\lambda+\gamma)\Delta}\mathbb{E}W(y(t_s+l\Delta))$. When $l=0$, based on system (\ref{15}) and Assumption 7, we have
\begin{eqnarray}
	&&\mathbb{E}W(y(t))\nonumber\\
	&&\leq e^{(\lambda+\gamma)\Delta}\mathbb{E}W(y(t_s))\nonumber\\
	&&=e^{(\lambda+\gamma)\Delta+\gamma t_s}\mathbb{E}[y^\mathrm{T}(t_s^--\tau_s)(e^{-\alpha\Gamma(t_s)}\tilde{K}\tilde{H}\Phi\Theta)^\mathrm{T}P\nonumber\\
	&&\quad\times (e^{-\alpha\Gamma(t_s)}\tilde{K}\tilde{H}\Phi\Theta)y(t_s^--\tau_s)]\nonumber\\
	&&\leq \tilde{\lambda}_1e^{(\lambda+\gamma)\Delta-2\alpha\beta t_s+\gamma\tau_s}\mathbb{E}W(y(t_s^--\tau_s)).\label{20}
\end{eqnarray}

When $l\in(0,l_s)$, according to PETM (\ref{07}) and  Assumption 7, we get
\begin{eqnarray}
	&&\mathbb{E}W(y(t))\nonumber\\
	&&\leq e^{(\lambda+\gamma)\Delta}\mathbb{E}W(y(t_s+l\Delta))\nonumber\\
	&&\leq\psi_2e^{(\lambda+\gamma)\Delta}\mathbb{E}W(y(t_s))\nonumber\\
	&&\leq\psi_2\tilde{\lambda}_1e^{(\lambda+\gamma)\Delta-2\alpha\beta t_s+\gamma\tau_s}\mathbb{E}W(y(t_s^--\tau_s)).\label{21}
\end{eqnarray}

As a matter of fact, (\ref{21}) is not only tenable for $t\in[t_s+l\Delta, t_s+(l+1)\Delta)$ but also suitable for $t\in[t_s+l_s\Delta, t_{s+1})$.

According to  $\psi_2\geq1$, $t_s\geq\Delta$ and Assumption 8, it follows from (\ref{20}) and (\ref{21}) that for $t\in[t_s, t_{s+1})$,
\begin{eqnarray}
	&&\mathbb{E}W(y(t))\nonumber\\
	&&\leq\psi_2\tilde{\lambda}_1e^{(\lambda+\gamma)\Delta-2\alpha\beta t_s+\gamma\tau_s}\mathbb{E}W(y(t_s^--\tau_s))\nonumber\\
	&&\leq\psi_2\tilde{\lambda}_1e^{(\lambda+2\gamma-2\alpha\beta) \Delta}\mathbb{E}W(y(t_s^--\tau_s)).\label{22}
\end{eqnarray}

Based on Assumption 9 and $\gamma\in(0,\bar{\gamma}^{'}]$, (\ref{22}) implies $\mathbb{E}W(y(t))\leq\mathbb{E}W(y(t_s^--\tau_s))$. Similar to (\ref{14}), for $t\geq t_0$, we have $\mathbb{E}W(y(t))\leq\mathbb{E}W(y(t_0^--\tau_0))$.

Under Assumption 8, we know $t_0^--\tau_0\in[t_0-\Delta,t_0)$. Then, applying (\ref{0014}), we obtain $\mathbb{E}|y(t)|^2\leq Me^{-\gamma t}\mathbb{E}|y(0)|^2.$ Therefore, this completes the proof of Theorem 2.
{\hfill $\Box$}

\begin{remark}
	Different from MASs investigated in \cite{HZWR23, LXLF22, KRSSYD24, YLCLY25, HZH25, LSL24}, we consider a fully heterogeneous stochastic MAS with a time-varying topology. In this paper, a virtual state space is constructed to unify the dimensions of the agents and ensure information interaction between them. Furthermore, based on the energy consumption in the time-varying topology, this paper proposes a novel PETM that conserves bandwidth resources. Finally, by designing PETICs with/without actuation delays, the consensus of the system is achieved.
\end{remark}

\section{Numerical simulation}

\begin{figure}[H]
	\begin{center}
		\includegraphics[height=5cm]{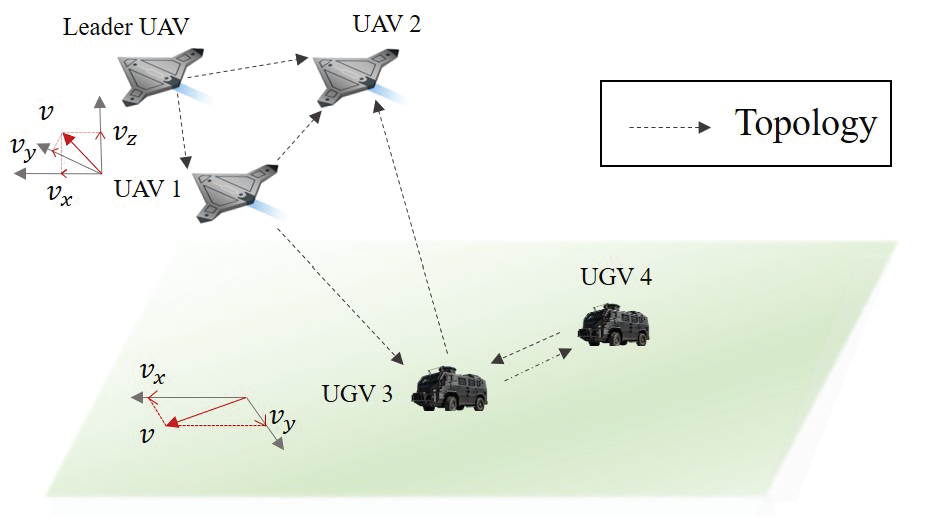}
		\caption{UAVs-UGVs system.}
		\label{topology}
	\end{center}
\end{figure}

In this section, we consider a leader-following system, where the leader is an UAV  and the followers consist of two UAVs and two UGVs. The topological diagram is shown as Fig. \ref{topology}. Based on the work in \cite{SWW25}, we incorporate stochastic noise into the system. Then,  the dynamic leader model is described as
\begin{eqnarray*}
	&&d\eta_0(t)=v_0(t)dt+d_0v_0(t)d\omega(t),\nonumber\\
	&&dv_0(t)=[c_{\eta0}\eta_0(t)+c_{v0}v_0(t)]dt\\
	&&\qquad\quad\quad+[d_{\eta0}\eta_0(t)+d_{v0}v_0(t)]d\omega(t),~t\geq 0,
\end{eqnarray*}
where $\eta_0(t)=[\eta_{0x}^\mathrm{T}(t),\eta_{0y}^\mathrm{T}(t),\eta_{0z}^\mathrm{T}(t)]^\mathrm{T}\in\mathbb{R}^3$ and $v_0(t)=[v_{0x}^\mathrm{T}(t),v_{0y}^\mathrm{T}(t),v_{0z}^\mathrm{T}(t)]^\mathrm{T}\in\mathbb{R}^3$ are the position and the velocity of leader UAV. Let $\bar{x}_0(t)=[\eta_0^\mathrm{T}(t),v_0^\mathrm{T}(t)]^\mathrm{T}$, the leader system is established as
\begin{eqnarray*}
	&&d\bar{x}_0(t)=\mathbb{C}_0\bar{x}_0(t)dt+\mathbb{D}_0\bar{x}_0(t)d\omega(t),~t\geq 0,
\end{eqnarray*}
where $\mathbb{C}_0\in\mathbb{R}^{6\times 6}$ and  $\mathbb{D}_0\in\mathbb{R}^{6\times 6}$.  

In what follows, a dynamic follower model is established as
\begin{eqnarray*}
	&&d\eta_i(t)=v_i(t)dt+d_iv_i(t)d\omega(t),\nonumber\\
	&&dv_i(t)=[c_{\eta i}\eta_i(t)+c_{vi}v_i(t)+f_i(\eta_i(t))]dt\\
	&&\qquad\quad\quad+[d_{\eta i}\eta_i(t)+d_{vi}v_i(t)]d\omega(t), ~i\in \mathcal{N},~t\geq0,
\end{eqnarray*}
where $\mathcal{N}=\{1,2,3,4\}$, $\eta_i(t)=[\eta_{ix}^\mathrm{T}(t),\eta_{iy}^\mathrm{T}(t),\eta_{iz}^\mathrm{T}(t)]^\mathrm{T}\in\mathbb{R}^3$ and $v_i(t)=[v_{ix}^\mathrm{T}(t),v_{iy}^\mathrm{T}(t),v_{iz}^\mathrm{T}(t)]^\mathrm{T}\in\mathbb{R}^3$, $i=1,2$, are the position and the velocity of  UAVs. $\eta_i(t)=[\eta_{ix}^\mathrm{T}(t),\eta_{iy}^\mathrm{T}(t)]^\mathrm{T}\in\mathbb{R}^2$ and $v_i(t)=[v_{ix}^\mathrm{T}(t),v_{iy}^\mathrm{T}(t)]^\mathrm{T}\in\mathbb{R}^2$, $i=3,4$, are the position and the velocity of UGVs.
Let $\bar{x}_i(t)=[\eta_i^\mathrm{T}(t),v_i^\mathrm{T}(t)]^\mathrm{T}$, we constructed the following systems:
\begin{eqnarray*}
	&&d\bar{x}_i(t)=[\mathbb{C}_i\bar{x}_i(t)+\bar{f}_i(\bar{x}_i(t))+\bar{u}_i(t)]dt+\mathbb{D}_i\bar{x}_i(t)d\omega(t),\\
	&&\quad\quad\quad\quad\quad\quad\quad\quad\quad\quad\quad\quad\quad\quad\qquad\quad\:\: i\in \mathcal{N},~t\geq0,
\end{eqnarray*}
where $\bar{f}_i(\bar{x}_i(t))=[\mathbf{0},f_i(\bar{x}_i(t))]^\mathrm{T}$, $\bar{u}_i(t)=[\mathbf{0},u_i(t)]^\mathrm{T}$,
$\mathbb{C}_i\in\mathbb{R}^{6\times 6}$ and  $\mathbb{D}_i\in\mathbb{R}^{6\times 6}$.

Let $\bar{z}_i(t)=\bar{x}_i(t)-\Xi_i\bar{x}_0(t)-\bar{x}_i^*$ and $\bar{y}_i(t)=\Phi_i\bar{z}_i(t)$, where $\bar{x}_i^*$ denotes the desired relative position and $\bar{x}_i^*=[\bar{x}^{\mathrm{T}}_{x i},\bar{x}^{\mathrm{T}}_{y i},\bar{x}^{\mathrm{T}}_{z i},\mathbf{0},\mathbf{0},\mathbf{0}]^\mathrm{T}$. It is easy to obtain that
\begin{align*}
	&d\bar{y}_i(t)=[\Phi_i\mathbb{C}_i\bar{y}_i(t)+\Phi_i\bar{f}_i(t)+\Phi_i\bar{u}_i(t)]dt\nonumber\\
	&\quad\quad\quad\quad+[\Phi_i\mathbb{D}_i\Theta_i\bar{y}_i(t)+\Phi_i\mathbb{D}_i\bar{x}_i^*]d\omega(t),~i\in \mathcal{N},~t\geq0,
\end{align*}
where $\bar{f}_i(t)=\bar{f}_i(\bar{z}_i(t))+\mathbb{C}_i\bar{x}_i^*$  and $|\Phi_i\mathbb{D}_i\Theta_i\bar{y}_i(t)+\Phi_i\mathbb{D}_i\bar{x}_i^*|\leq L_{Di}|\bar{y}_i(t)|$.

Then, the dynamic parameters matrix of UAVs-UGVs system can be described as:
\begin{equation*}
	\begin{aligned}
		&\mathbb{C}_0=\ diag\left\{
		\setlength{\arraycolsep}{1.8pt}
		\begin{bmatrix}0&&1\\0&&0\end{bmatrix},
		\begin{bmatrix}0&&1\\0&&0\end{bmatrix},
		\begin{bmatrix}0&&1\\-0.05&&0.05\end{bmatrix}
		\right\},\\[6pt]
		&\mathbb{D}_0=\ diag\left\{
		\setlength{\arraycolsep}{1.8pt}
		\begin{bmatrix}0&&0.01\\0&&0\end{bmatrix},
		\begin{bmatrix}0&&0.01\\0&&0\end{bmatrix},
		\begin{bmatrix}0&&0.01\\0&&0\end{bmatrix}
		\right\},\\[6pt]
		&H=\setlength{\arraycolsep}{3pt}		
		\begin{bmatrix}
			1&   0 &  0 &  0\\
			1  &  -1   &     1   &     0\\
			1  & 0   &  -2 & 1\\
			0  &     0      &  1 &   -1
		\end{bmatrix}.
	\end{aligned}
\end{equation*}

For $i=1,~2$,
\begin{equation*}
	\begin{aligned}
		&\mathbb{C}_i=\ diag\left\{
		\setlength{\arraycolsep}{1.8pt}
		\begin{bmatrix}0&&1\\-0.9&&0.9\end{bmatrix},
		\begin{bmatrix}0&&1\\0.6&&0.6\end{bmatrix},
		\begin{bmatrix}0&&1\\-0.01&&0.01\end{bmatrix}
		\right\},\\[6pt]
		&\mathbb{D}_i=\ diag\left\{
		\setlength{\arraycolsep}{1.8pt}
		\begin{bmatrix}0&&1.5\\0.6&&0.6\end{bmatrix},
		\begin{bmatrix}0&&1.9\\1.6&&1.6\end{bmatrix},
		\begin{bmatrix}0&&1.7\\0.6&&0.6\end{bmatrix}
		\right\},
	\end{aligned}
\end{equation*}
$\Xi_i=\Phi_i=\Theta_i=I_6$, $\bar{f}_i(t)=[0,0,0,0.5sin(0.3\eta_{ix}(t)),0.7\\\times sin(0.5\eta_{iy}(t)),0.5sin(0.9\eta_{iz}(t))]^\mathrm{T}$, $\bar{x}_1^*=[3,3,0,0,0,0]^\mathrm{T}$ and $\bar{x}_2^*=[0,3,0,0,0,0]^\mathrm{T}$. Based on \cite{BKAT22, YCY20}, we normalize the parameters to obtain the energy consumption $\tau_i(t)=\tau+0.005t$.

For $i=3,~4$,
\begin{equation*}
	\begin{aligned}
	    &\mathbb{C}_i=diag\left\{\setlength{\arraycolsep}{1.8pt}\begin{bmatrix}0&&1\\-0.8&&0.8\end{bmatrix},\begin{bmatrix}0&&1\\-0.6&&0.6\end{bmatrix}\right\},\\
		&\mathbb{D}_i=diag\left\{\setlength{\arraycolsep}{1.8pt}\begin{bmatrix}0&&1.5\\1.2&&1.2\end{bmatrix},\begin{bmatrix}0&&1.9\\0.5&&0.5\end{bmatrix}\right\},\\
		&\Xi_i=\Theta_i=\begin{bmatrix}1&0&0&0&0&0\\0&1&0&0&0&0\\0&0&0&1&0&0\\0&0&0&0&1&0\end{bmatrix},~\Phi_i=\begin{bmatrix}1&0&0&0\\0&1&0&0\\0&0&0&0\\0&0&1&0\\0&0&0&1\\0&0&0&0\end{bmatrix},
	\end{aligned}
\end{equation*}
$\bar{f}_i(t)=[0,0,0.5sin(0.3\eta_{ix}(t)),0.7sin(0.5\eta_{iy}(t))]^\mathrm{T}$, $\bar{x}_3^*=[2,2,0,0]^\mathrm{T}$, $\bar{x}_4^*=[0,0,0,0]^\mathrm{T}$ and $\tau_i(t)=\tau+0.0047t$.

\begin{remark}
	Different from \cite{SWW25}, we introduce stochastic noise into the system, which more accurately reflects the actual situation. In a stochastic operating environment, the agents we consider exhibit characteristics across multiple dimensions. Thus, we can achieve the coordinated operation of UAVs and UGVs. Such coordinated operations of UAVs and UGVs have also been applied and verified in \cite{LL24, XXTLS25}.  However, it should be noted that neither of these studies takes stochastic noise into account. Due to the unobservable nature of stochastic noise, achieving information transmission between adjacent heterogeneous agents becomes challenging. This difficulty also demonstrates the necessity of constructing the virtual state space proposed in this paper.
\end{remark}

The simulation is performed by MATLAB with step $= 0.003$. The initial values are given as follows:
\begin{eqnarray*}
	&&\eta_0(0)=[1.55,-1.5,10.01]^\mathrm{T},~v_0(t)=[2.3,-3.15,1.1]^\mathrm{T},\\
	&&\eta_1(0)=[7.11,-2.78,2.4]^\mathrm{T},~v_1(0)=[2,1.52,0.51]^\mathrm{T},\\
	&&\eta_2(0)=[10.12,7.63,2.6]^\mathrm{T},~v_2(0)=[-0.5,-1.02,1.3]^\mathrm{T},\\
	&&\eta_3(0)=[10.51,-7.7]^\mathrm{T},~v_3(0)=[2,-2.92]^\mathrm{T},\\
	&&\eta_4(0)=[5.12,5.63]^\mathrm{T},~v_4(0)=[1.5,-3.02]^\mathrm{T}.
\end{eqnarray*}

\begin{figure}[H]
	\begin{center}
		\includegraphics[height=5cm]{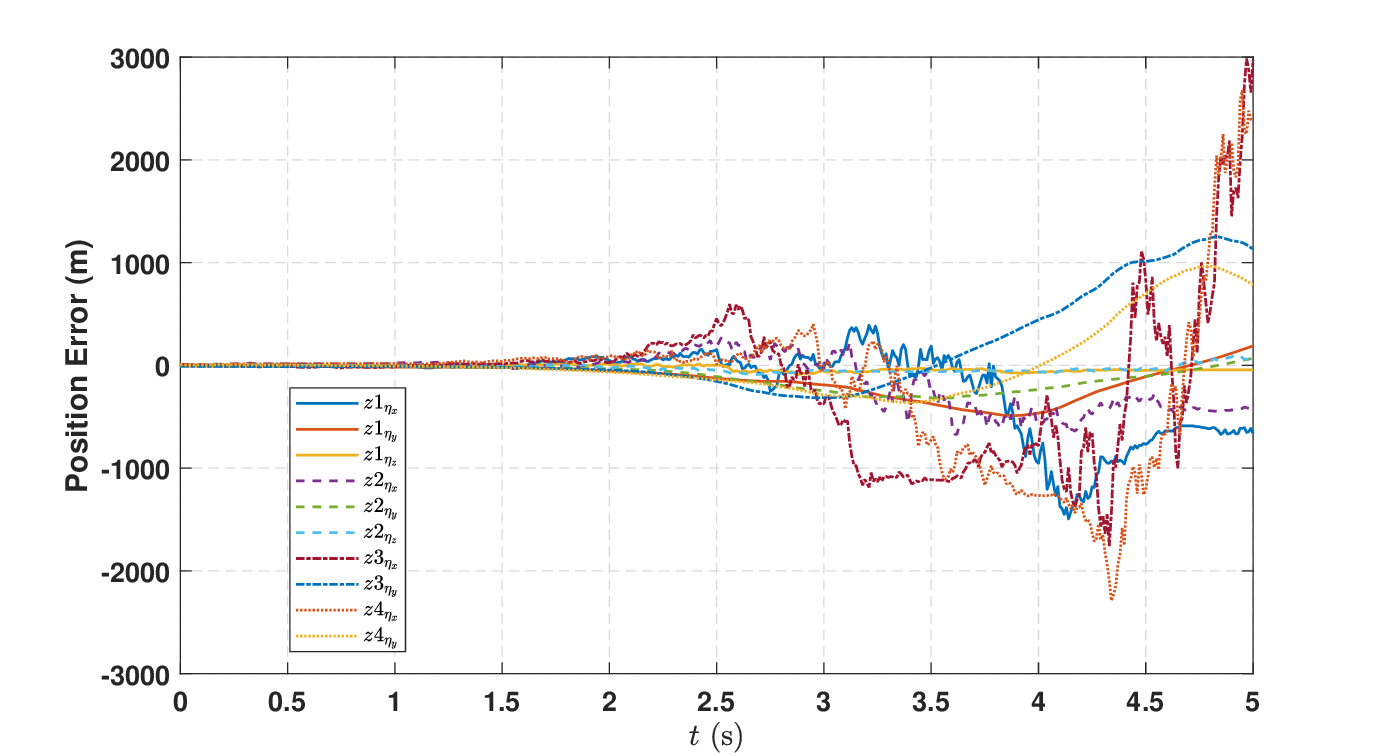}
		\caption{Position errors without control.}
		\label{eta0}
	\end{center}
\end{figure}

\begin{figure}[H]
	\begin{center}
		\includegraphics[height=5cm]{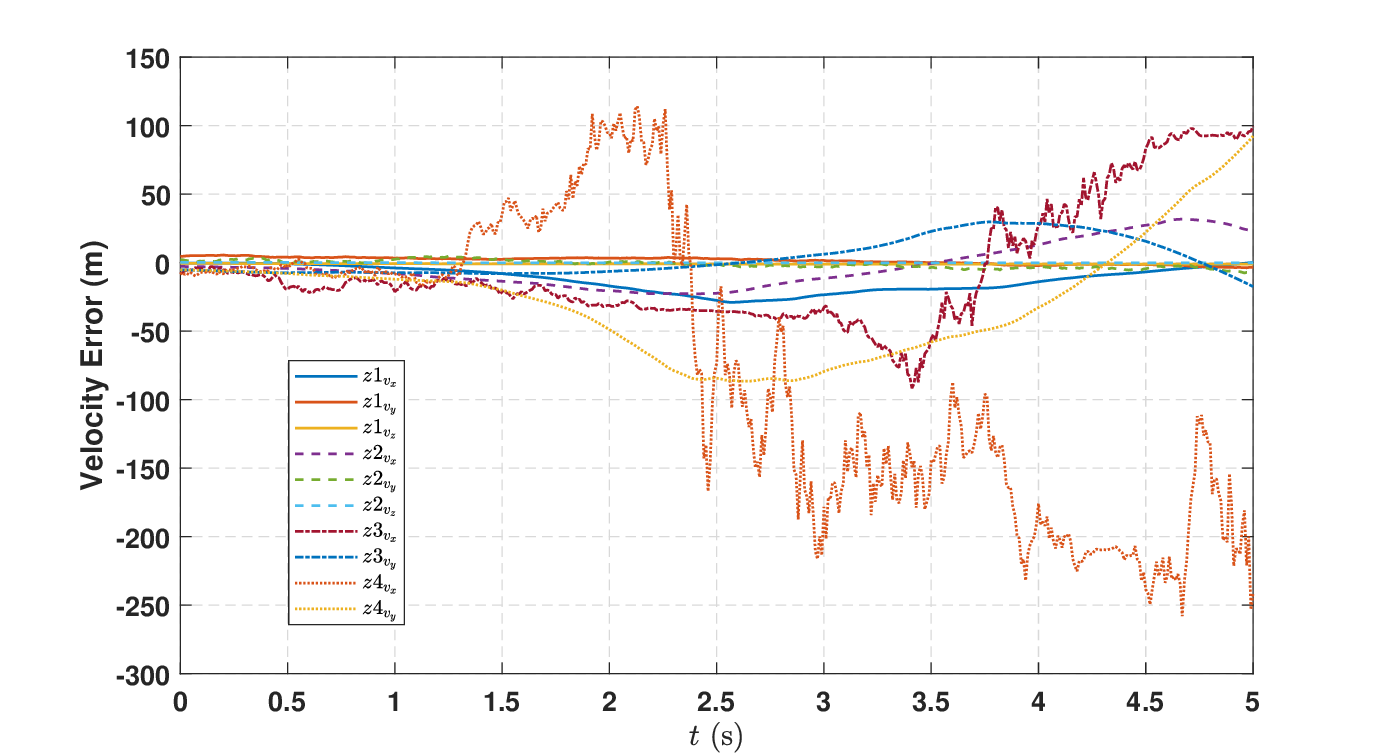}
		\caption{ Velocity errors without control.}
		\label{v0}
	\end{center}
\end{figure}

Then, Fig. \ref{eta0} and Fig. \ref{v0} illustrate the trajectories of UAVs-UGVs system without control. It can be seen that uncontrolled system is not stable.

First, we consider the case of control without actuator delays:
\begin{eqnarray}
	&&\bar{u}_i(t)=\Sigma_{s=0}^\infty k_i\bar{U}_i(t)\delta(t-t_s),\label{23}
\end{eqnarray}
where $\bar{U}_i(t)=e^{-\alpha\tau_i(t)}[\Sigma_{j= 1}^N(\bar{a}_{ij}\Theta_i(\bar{y}_i(t)-\bar{y}_j(t)))+\bar{b}_i\Theta_i\\\times\bar{y}_i(t)]$. 

According to system (\ref{06}), to verify the effectiveness of PETIC (\ref{23}), we choose parameters $k_1=-1.51$, $k_2=-2.38$, $k_3=-3.41$, $k_4=-1.10$, $\alpha=12.34$, $\tau=0.003$, $\psi_1=1.4$, $\psi_2=1.001$,  $\Delta=0.009$ and $P=0.95 I_{24}$. By calculation, we obtain $\lambda=10.99$, $\lambda_1=0.98$, which indicates that Assumption 6 is satisfied. Furthermore, we choose $\gamma=0.0019$ such that  $
\psi_2\lambda_1 e^{(\lambda+\gamma)\Delta-2\alpha\tau}\leq 1$. Therefore, all conditions in Theorem 1 are satisfied.

Then, we draw the following figures. Fig. \ref{1_y} depicts the trajectory of the virtual state $y(t)$ under PETIC (\ref{23}). Figs. \ref{1_z1}, \ref{1_z2}, \ref{1_z3} and \ref{1_z4} represent the state errors $z(t)$ of the four followers. Fig. \ref{1_intervals} plots the triggering intervals. Fig. \ref{1_3D} presents the dynamic operating trajectories of UAVs and UGVs during 20s. Based on these figures, it can be observed that the virtual states effectively reflect the movements of states with different dimensions. Meanwhile, it is evident that the state errors of UAVs and UGVs quickly achieve mean-square exponential stabilization under PETIC (\ref{23}).
Compared to the fixed-period control that requires $556$ controller updates, the number of updates is only $294$ by PETIC (\ref{23}), resulting in a $48\%$ reduction in data transmissions. Also, the Zeno behavior does not occur.

Second, we consider the case of control with actuator delays:
\begin{eqnarray}
	&&\tilde{u}_i(t)=\Sigma_{s=1}^\infty \tilde{k}_i\tilde{U}_i(t)\delta(t-t_s),\label{24}
\end{eqnarray}
where $\tilde{U}_i(t)=e^{-\alpha\tau_i(t)}[\Sigma_{j= 1}^N[\bar{a}_{ij}\Theta_i(y_i(t-\tau_s)-y_j(t-\tau_s))]+\bar{b}_i\Theta_iy_i(t-\tau_s)]$.

According to system (\ref{15}), we choose parameters $\tau_s=0.04$, $\tilde{k}_1=-0.71$, $\tilde{k}_2=0.55$, $\tilde{k}_3=0.08$, $\tilde{k}_4=-0.32$, $\alpha=2000$, $\tau=0$, $\beta=0.0047$, $\psi_1=1.2$, $\psi_2=1.1$, $\Delta=0.09$ and $P=0.5I_{24}$. By calculation, we obtain $\lambda=8.19$, $\tilde{\lambda}_1=1.16$, which indicates that Assumption 9 is satisfied. Based on Remark 13, we get $\bar{\gamma}^{'}=0.031$. Then, we choose $\gamma=0.03$. So, all conditions in Theorem 2 are satisfied.

 Fig. \ref{y} shows the behaviour of the virtual state $y(t)$ under PETIC (\ref{24}) with actuation delays. Figs. \ref{z1}, \ref{z2}, \ref{z3} and \ref{z4} plot the state errors $z(t)$. Fig. \ref{intervals} presents the time intervals of control with actuation delays. Fig. \ref{3D} depicts the dynamic operating trajectories of UAVs and UGVs. From these figures, it can be seen that UAVs and UGVs overcome the influence of actuator delays and achieve stability. What's more, in comparison with the fixed-period control that required about $55$ updates, the proposed PETIC (\ref{24}) achieves a $71\%$ reduction in updating frequency, requiring only $16$ updates.
 
 \begin{figure*}[!t]
 	\centering
 	\subfloat[Behaviours of virtual state $y(t)$ under PETIC (\ref{00}).]{\includegraphics[width=2.2in]{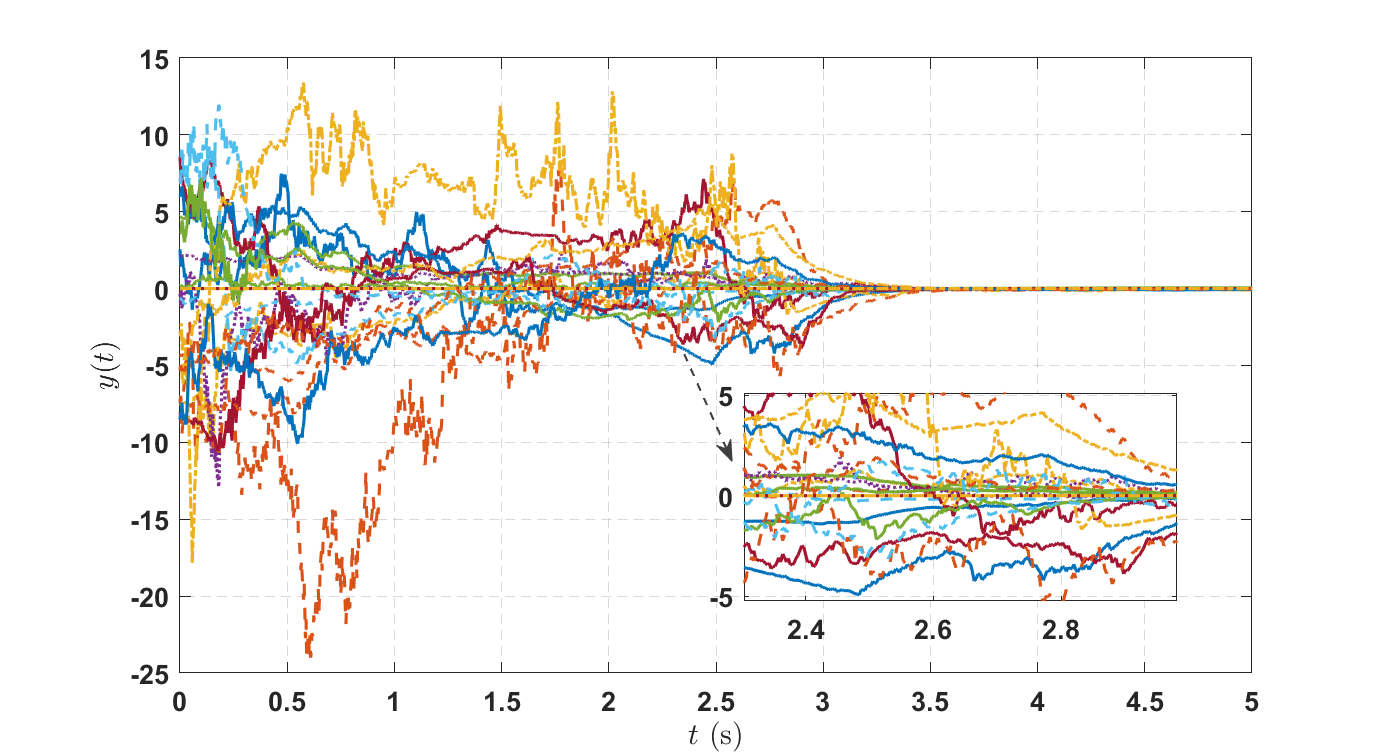}
 		\label{1_y}}
 	\hfil
 	\subfloat[Triggering intervals with PETM.]{\includegraphics[width=2.2in]{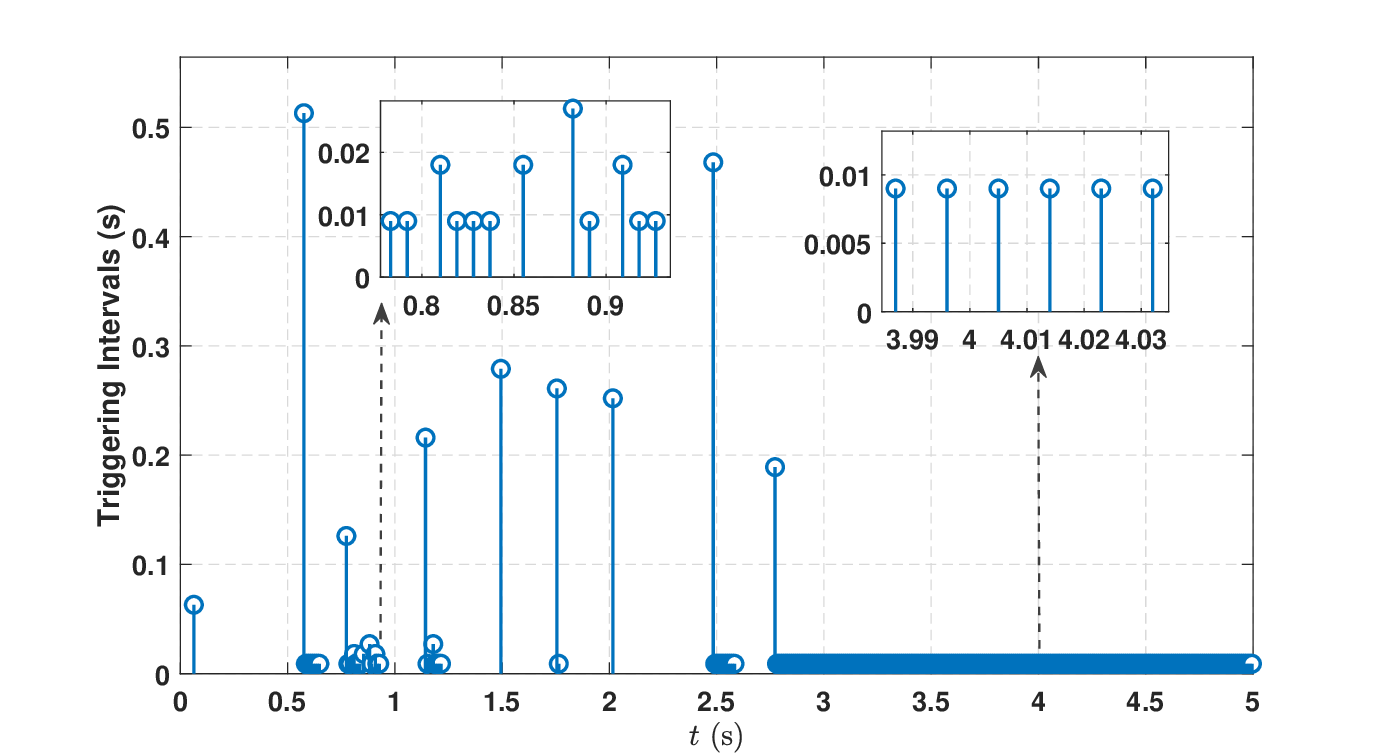}
 		\label{1_intervals}}
 	\hfil
 	\subfloat[The state error $z_1(t)$ of UAV 1.]{\includegraphics[width=2.2in]{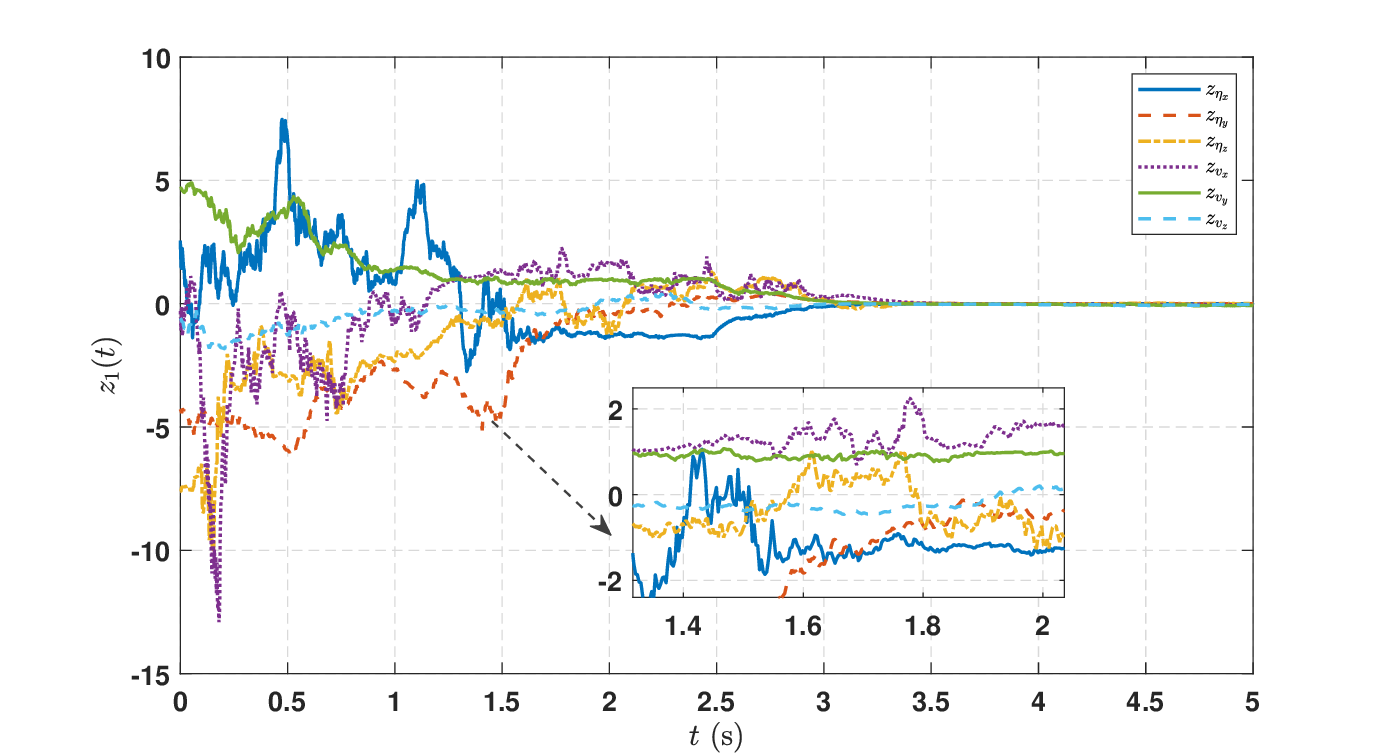}
 		\label{1_z1}}
 	\\
 	\subfloat[The state error $z_2(t)$ of UAV 2.]{\includegraphics[width=2.2in]{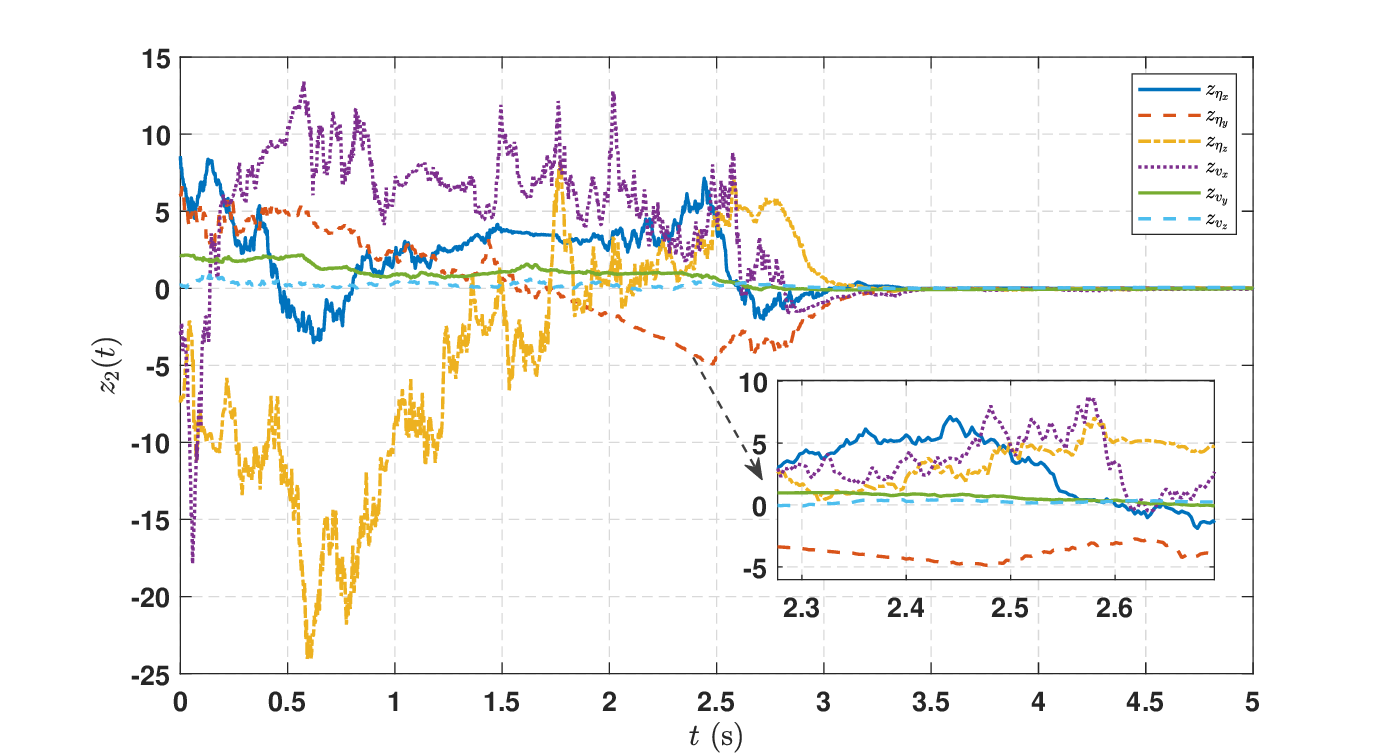}
 		\label{1_z2}}
 	\hfil
 	\subfloat[The state error $z_3(t)$ of UGV 3.]{\includegraphics[width=2.2in]{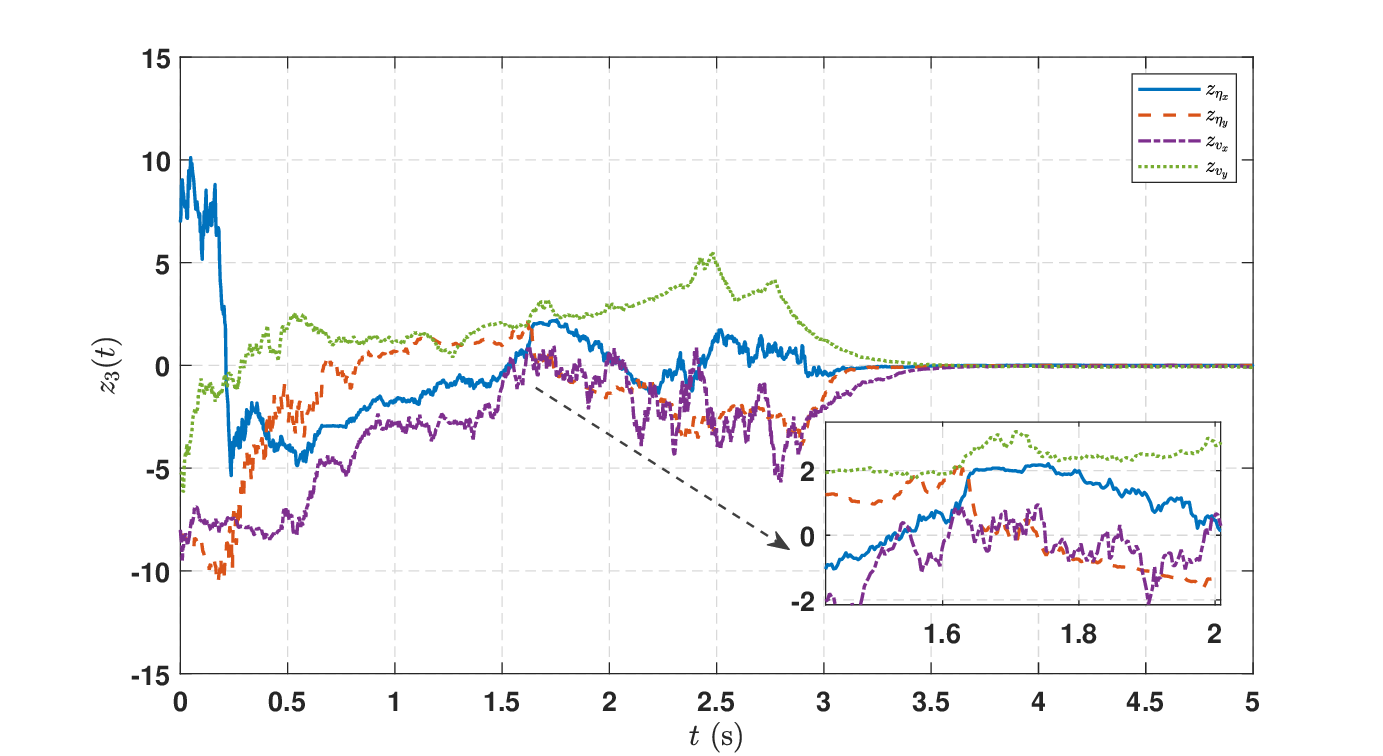}
 		\label{1_z3}}
 	\hfil
 	\subfloat[The state error $z_4(t)$ of UGV 4.]{\includegraphics[width=2.2in]{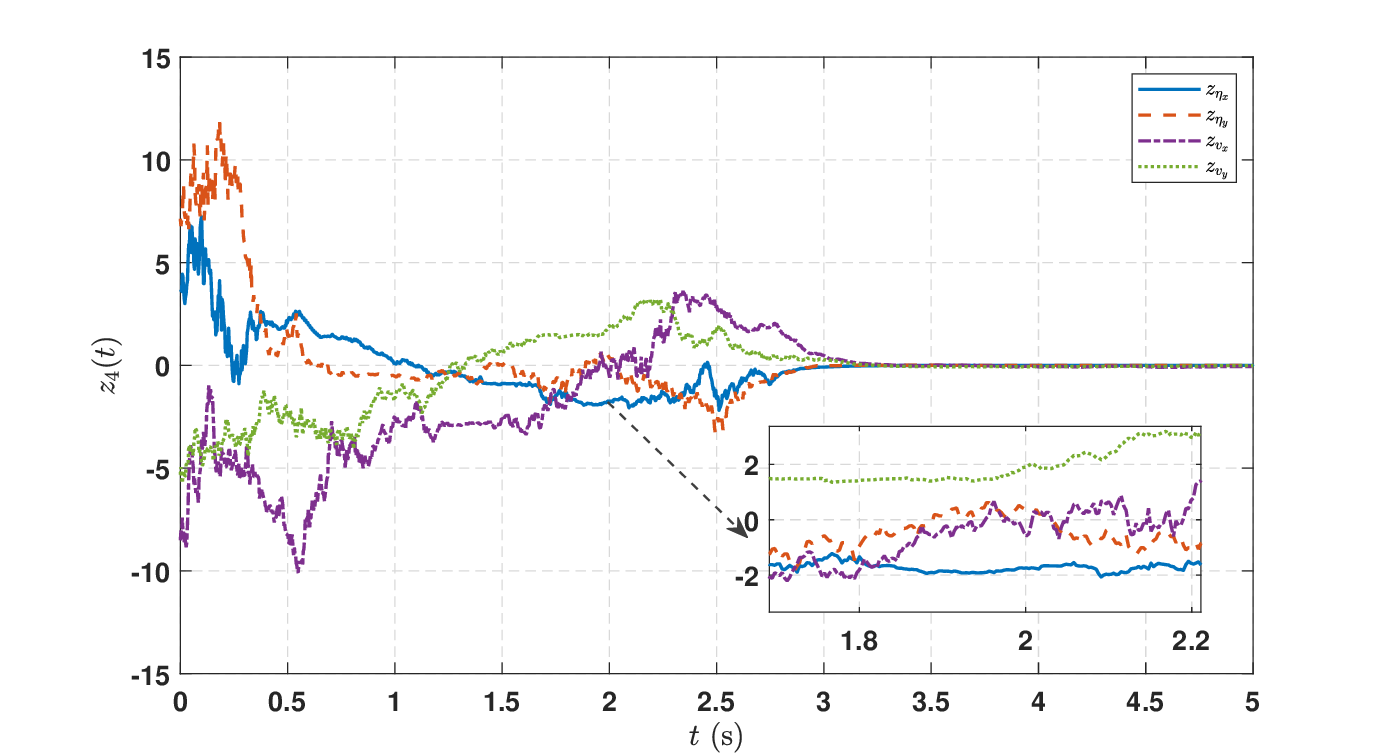}
 		\label{1_z4}}
 	\caption{Simulation results of UAVs and UGVs under PETIC (\ref{00}) without actuation delays.}
 	\label{1_simulation}
 \end{figure*}
 \begin{figure*}[!t]
 	\centering
 	\subfloat[Dynamic trajectories during 5s.]{\includegraphics[width=2.2in]{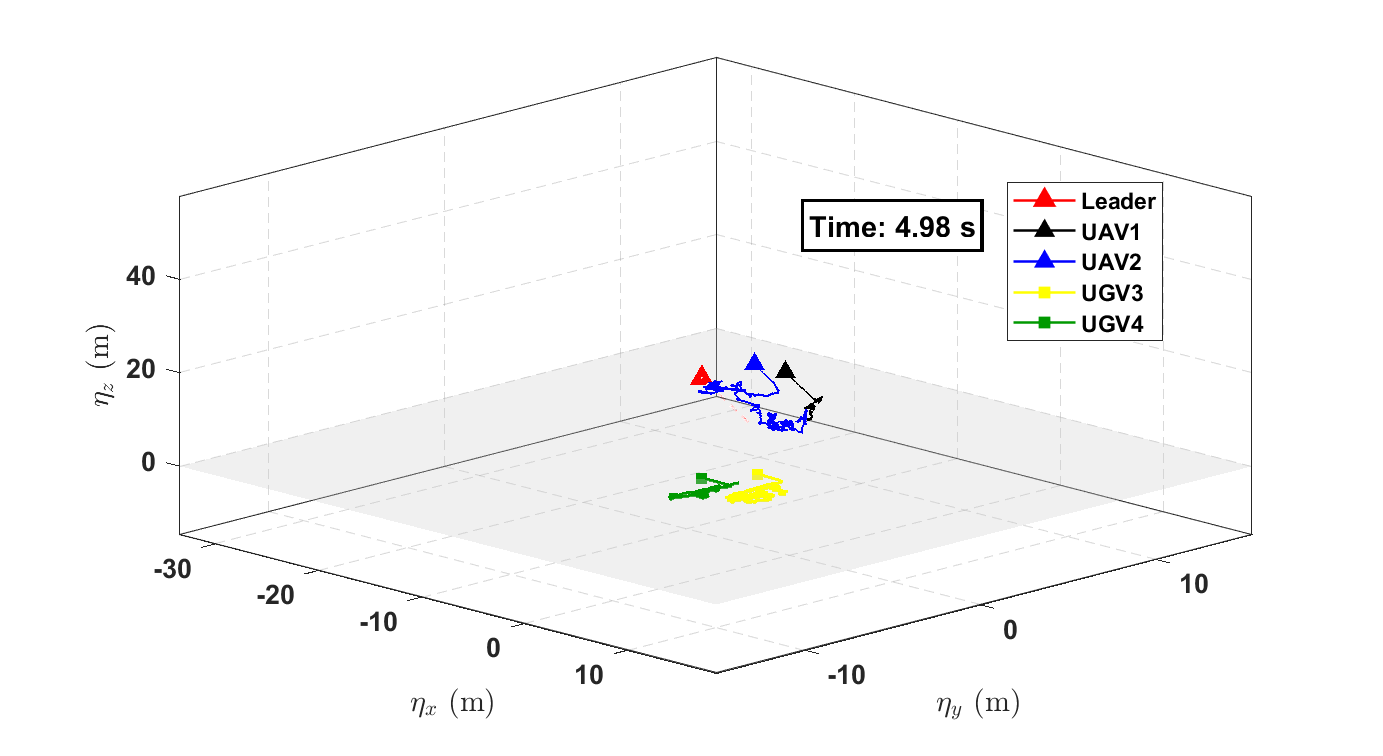}
 		\label{131}}
 	\hfil
 	\subfloat[Dynamic trajectories during 10s.]{\includegraphics[width=2.2in]{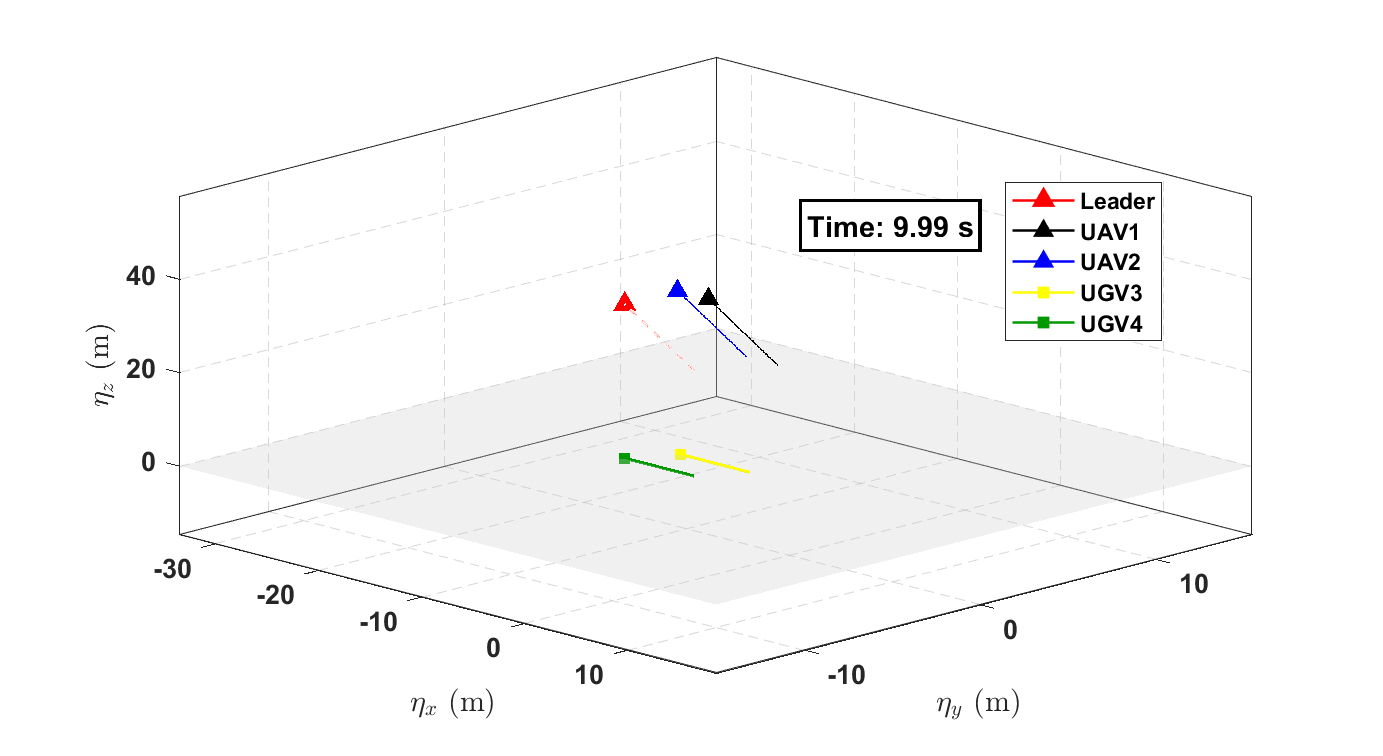}
 		\label{132}}
 	\hfil
 	\subfloat[Dynamic trajectories during 20s.]{\includegraphics[width=2.2in]{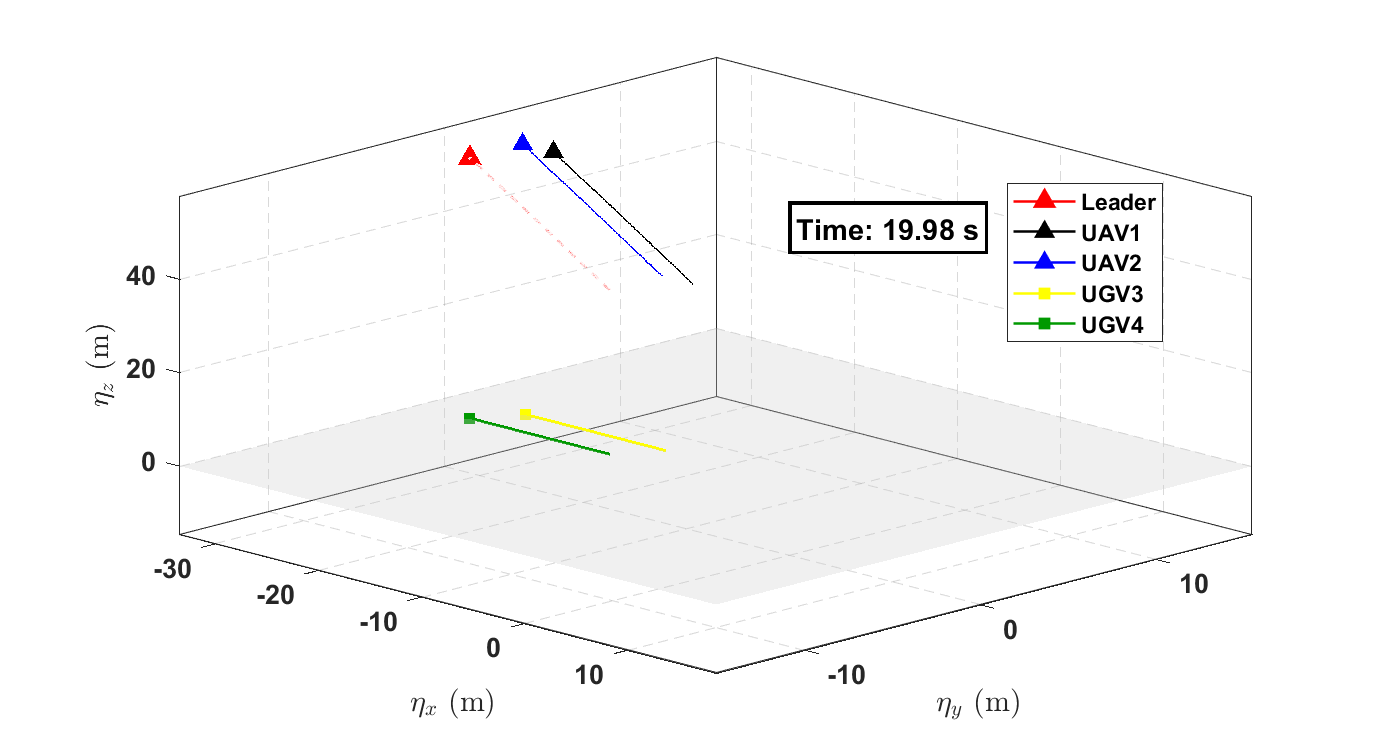}
 		\label{133}}
 	\caption{Dynamic trajectories of UAVs and UGVs under PETIC (\ref{00}) without actuation delays. }
 	\label{1_3D}
 \end{figure*}
\begin{figure*}[!t]
	\centering
	\subfloat[Behaviours of virtual state $y(t)$ under PETIC (\ref{000}).]{\includegraphics[width=2.2in]{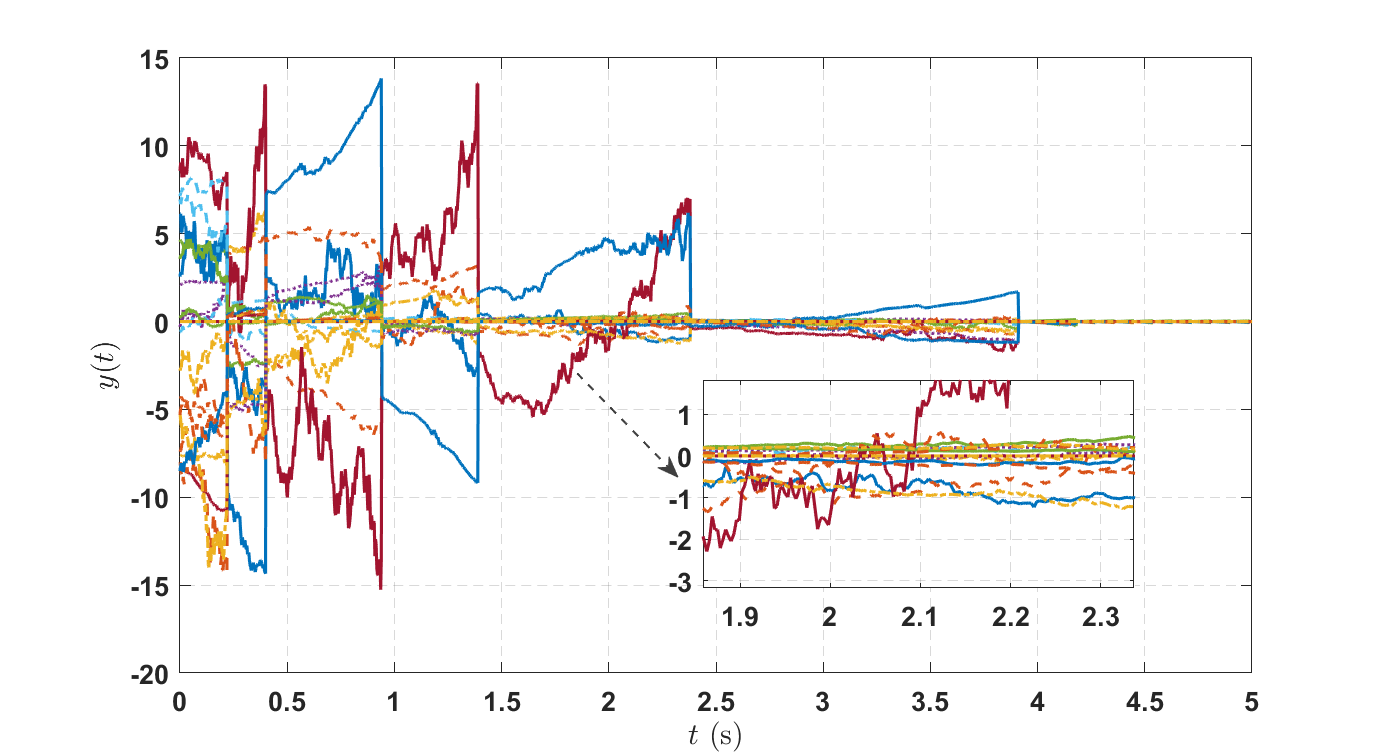}
		\label{y}}
	\hfil
	\subfloat[Control intervals with PETM.]{\includegraphics[width=2.2in]{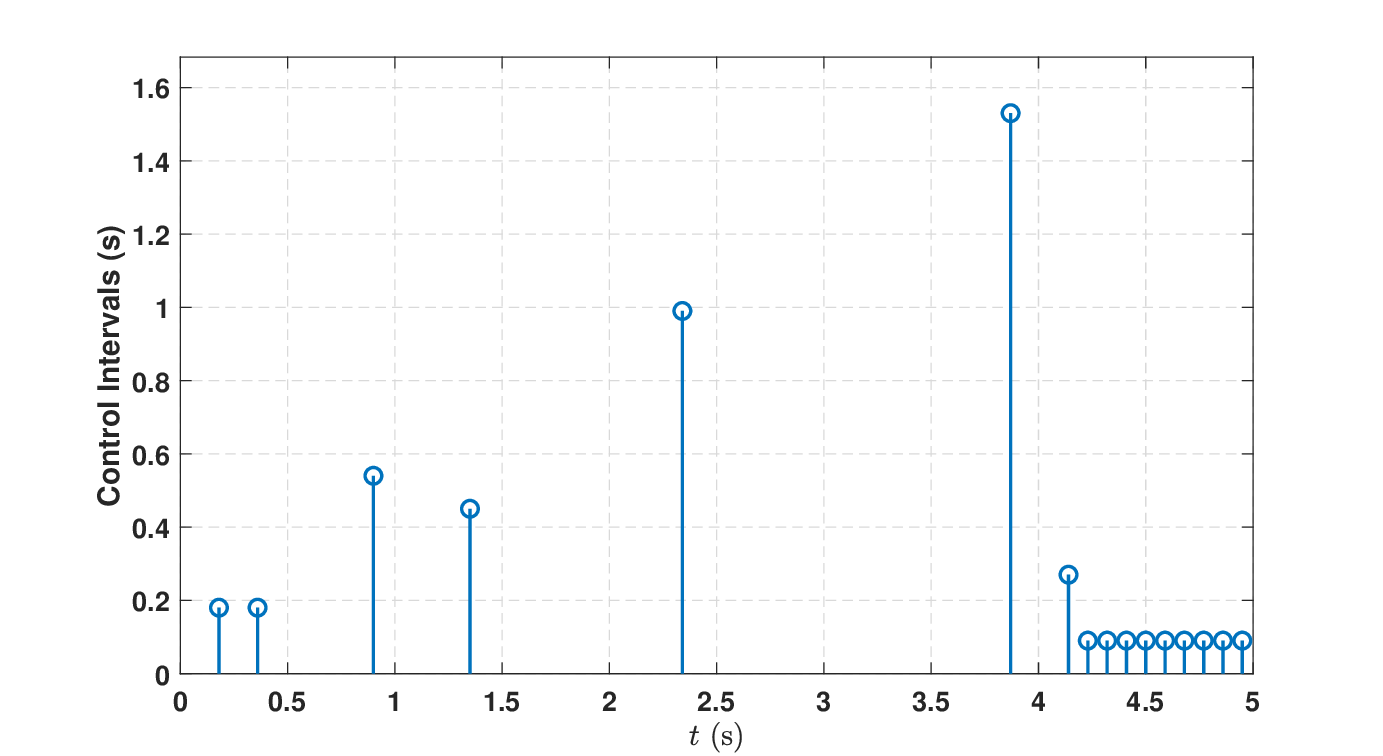}
		\label{intervals}}
	\hfil
	\subfloat[The state error $z_1(t)$ of UAV 1.]{\includegraphics[width=2.2in]{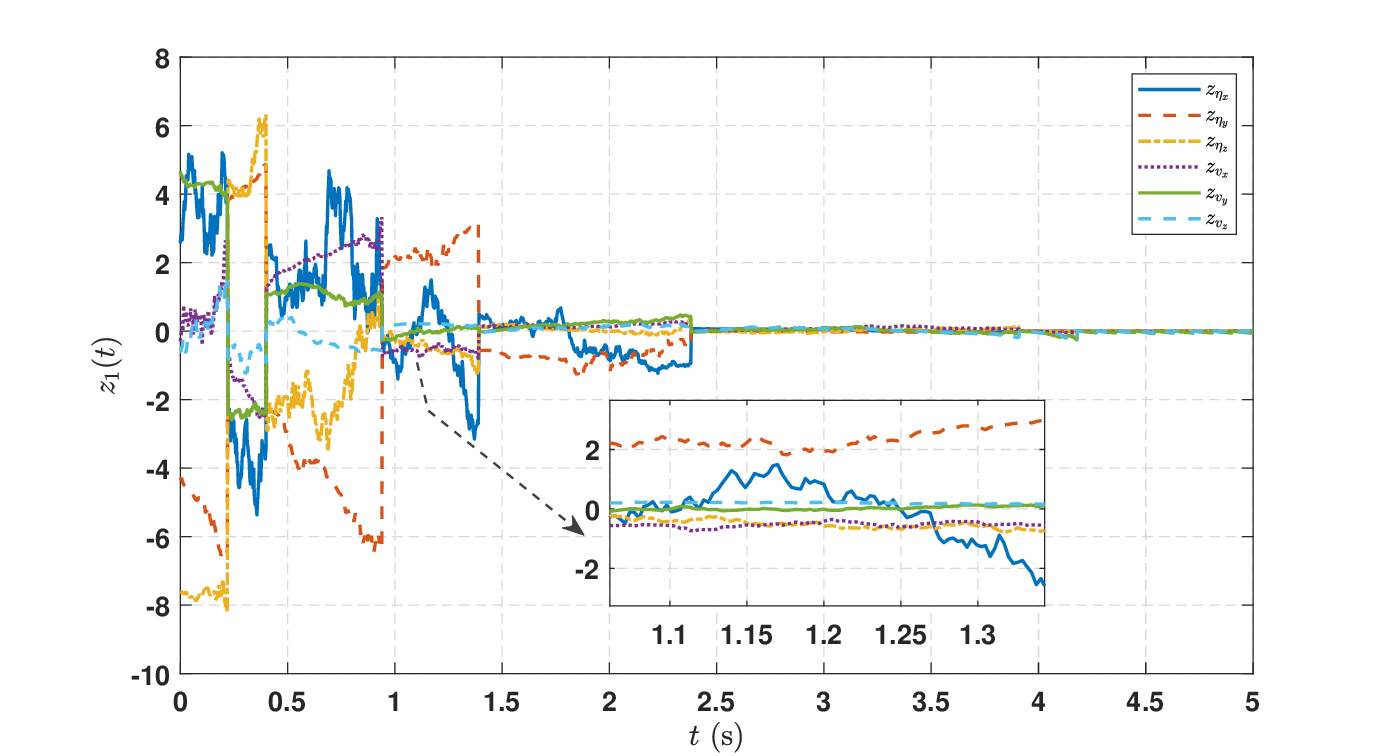}
		\label{z1}}
	\\
	\subfloat[The state error $z_2(t)$ of UAV 2.]{\includegraphics[width=2.2in]{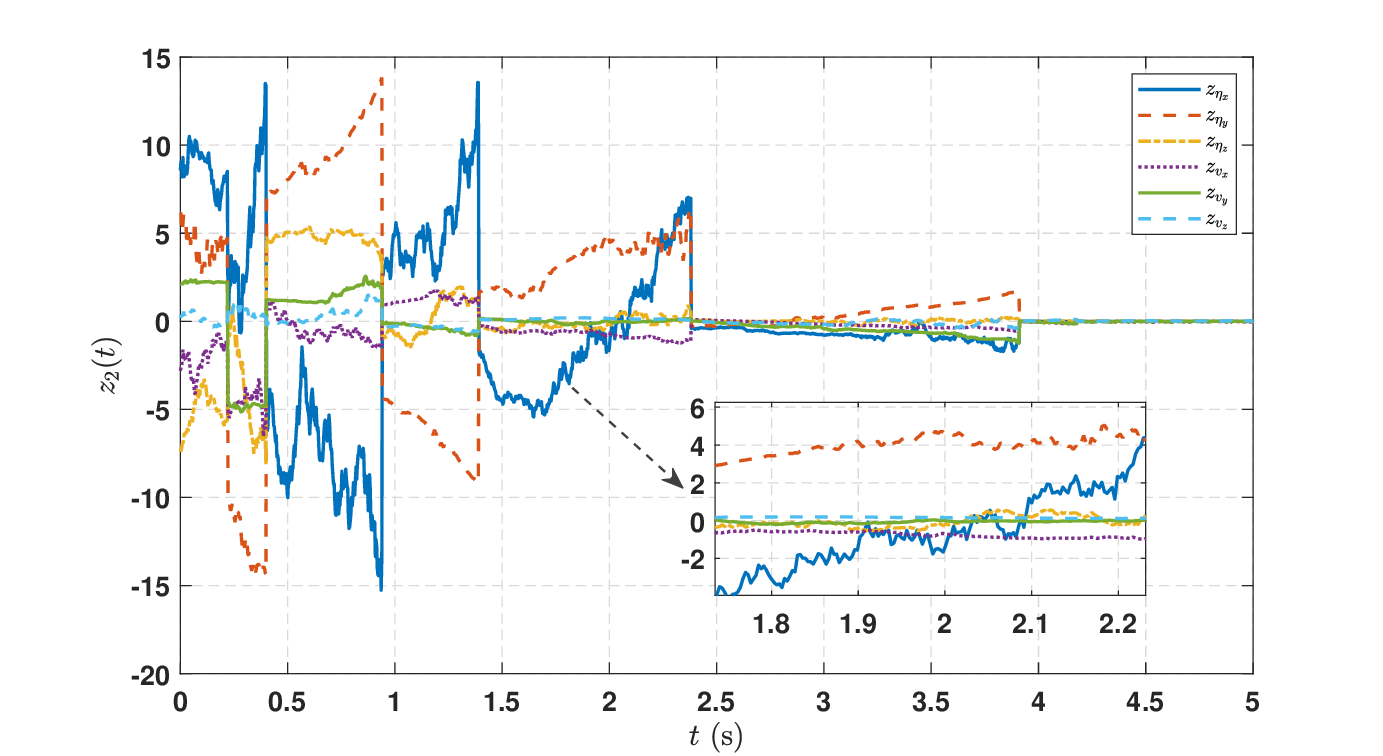}
		\label{z2}}
	\hfil
	\subfloat[The state error $z_3(t)$ of UGV 3.]{\includegraphics[width=2.2in]{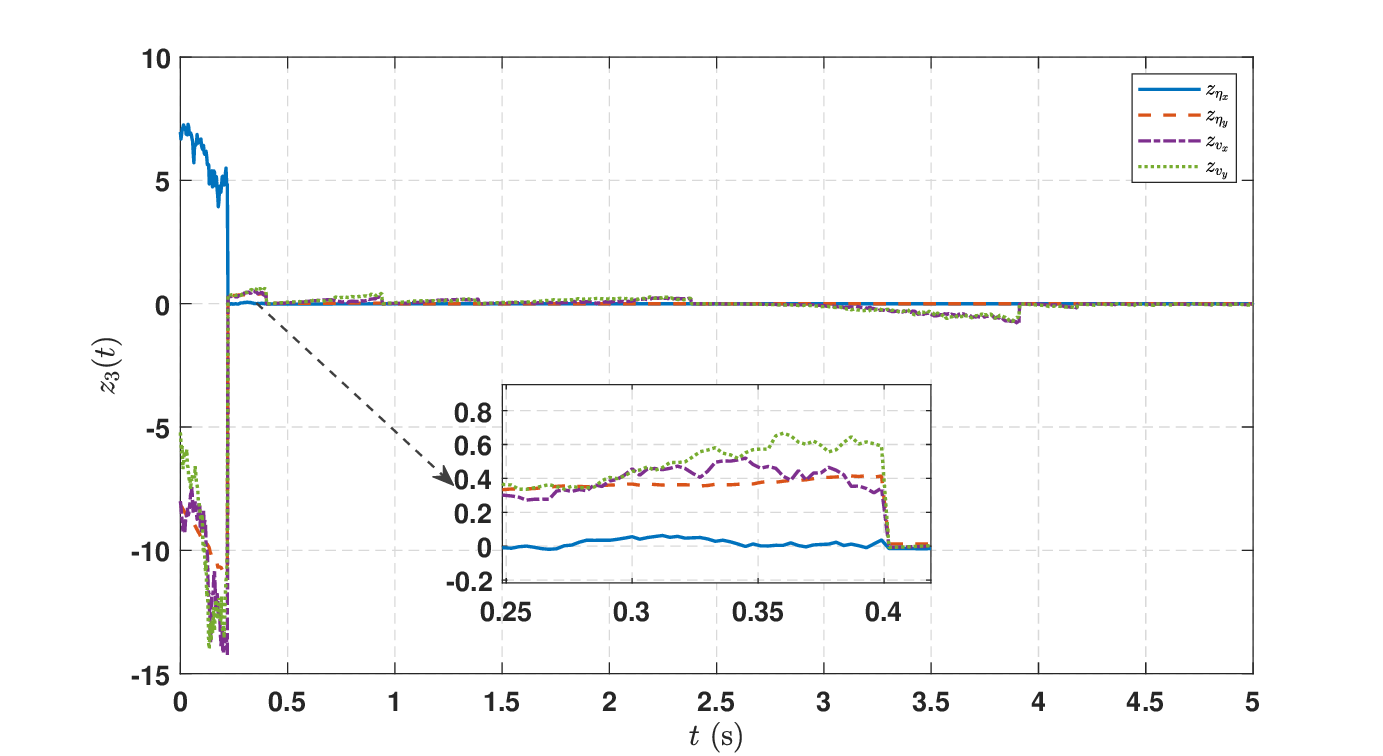}
		\label{z3}}
	\hfil
	\subfloat[The state error $z_4(t)$ of UGV 4.]{\includegraphics[width=2.2in]{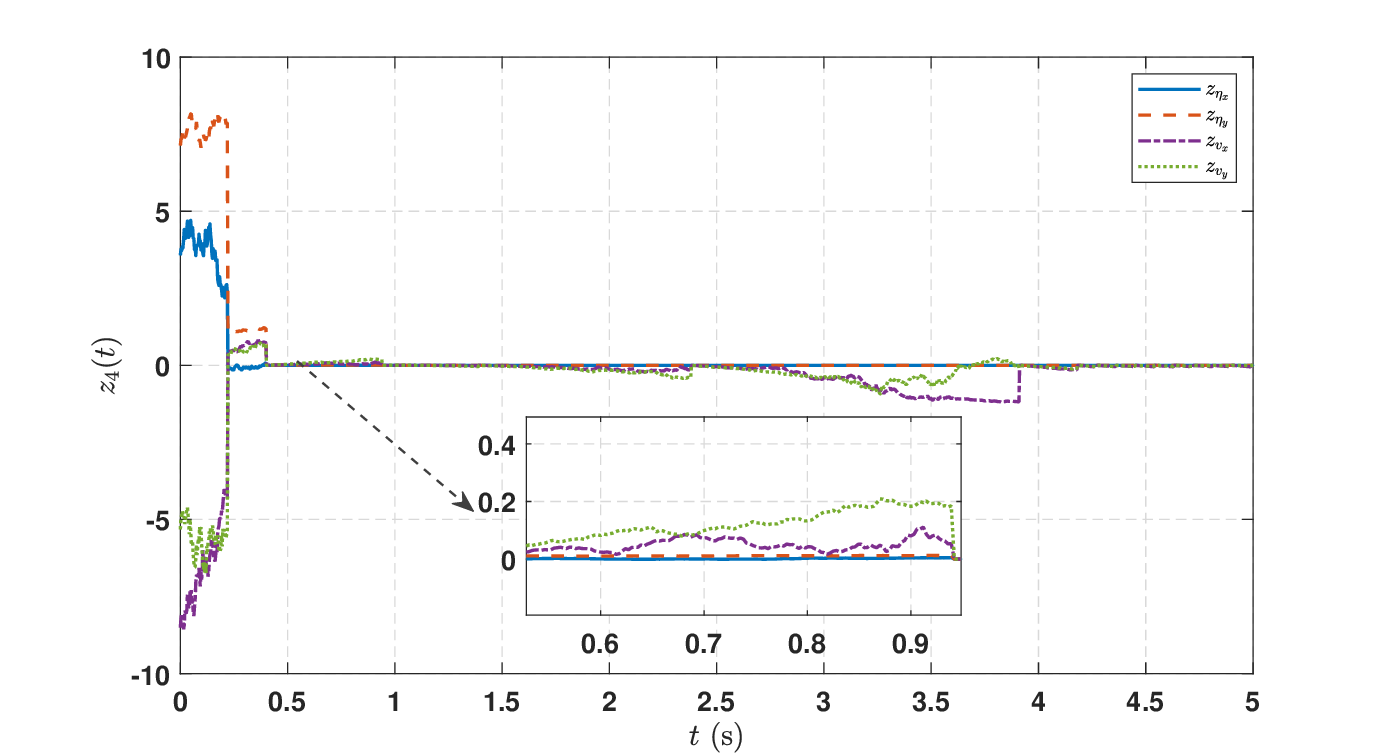}
		\label{z4}}
	\caption{Simulation results of UAVs and UGVs under PETIC (\ref{000}) with actuation delays. }
	\label{simulation}
\end{figure*}
\begin{figure*}[!t]
	\centering
	\subfloat[Dynamic trajectories during 5s.]{\includegraphics[width=2.2in]{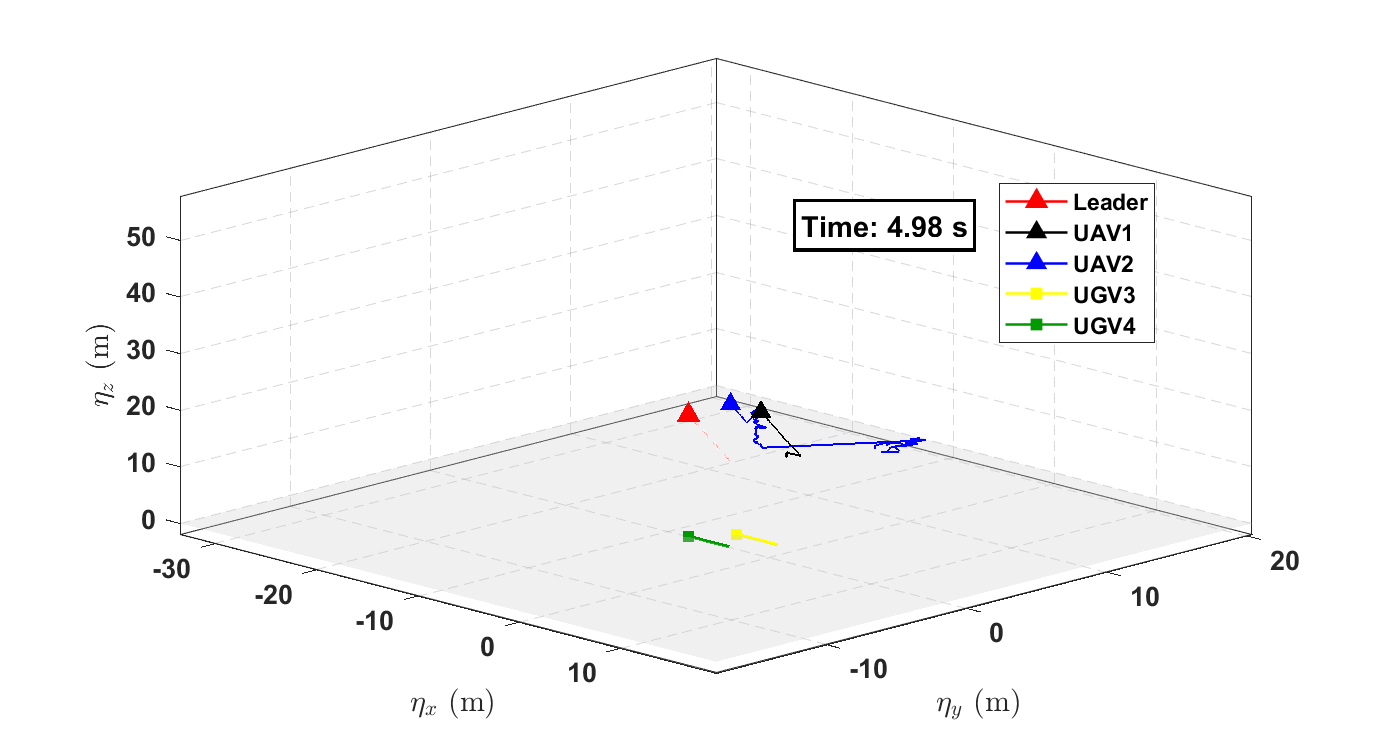}
		\label{31}}
	\hfil
	\subfloat[Dynamic trajectories during 10s.]{\includegraphics[width=2.2in]{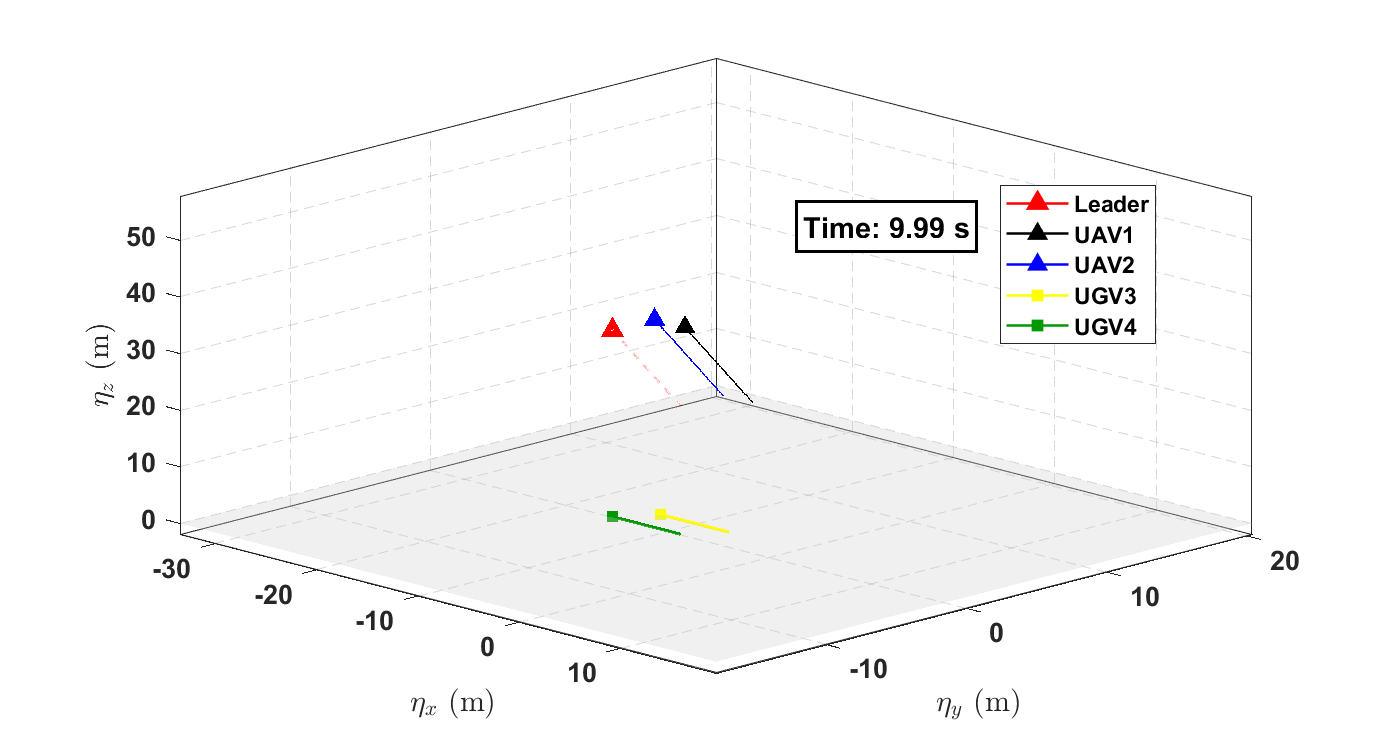}
		\label{32}}
	\hfil
	\subfloat[Dynamic trajectories during 20s.]{\includegraphics[width=2.2in]{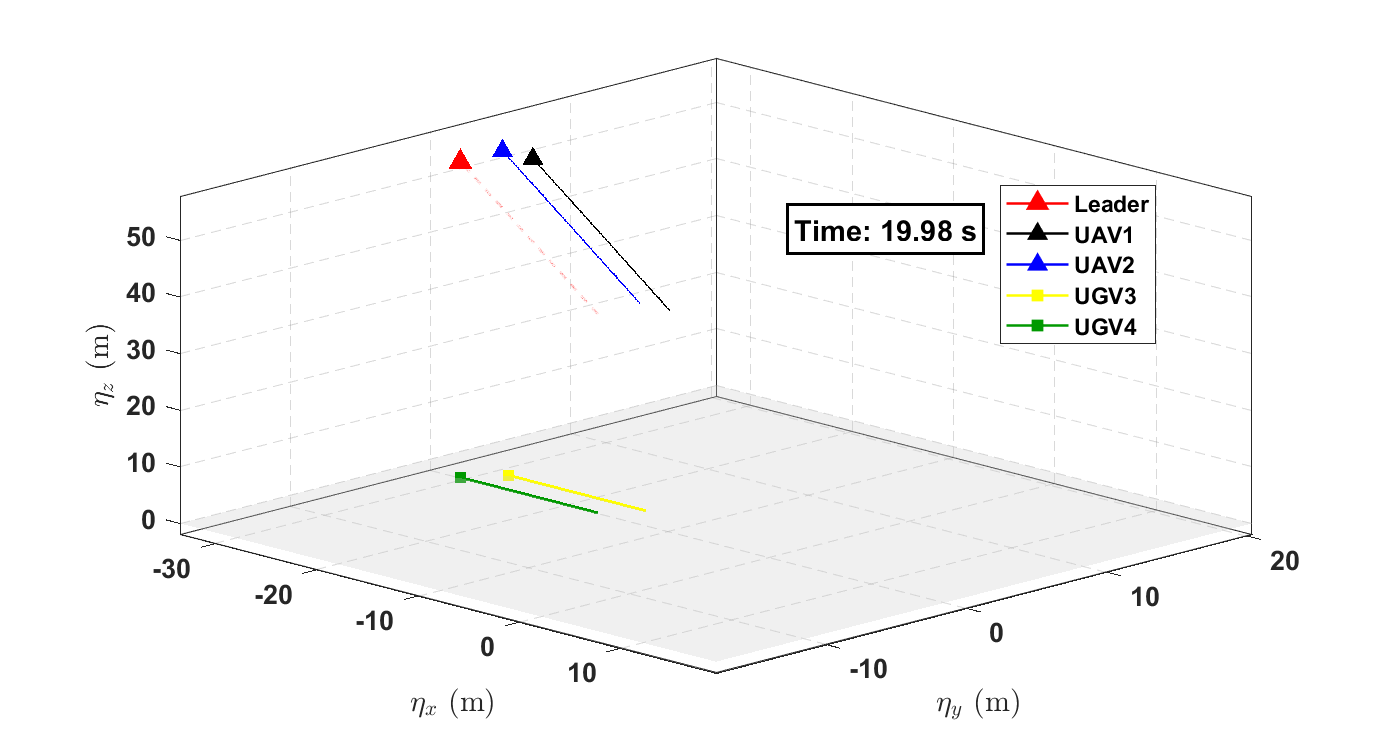}
		\label{33}}
	\caption{Dynamic trajectories of UAVs and UGVs under PETIC (\ref{000}) with actuation delays. }
	\label{3D}
\end{figure*}

\section{Conclusions and future research prospects}

Different from the traditional stochastic MASs, the time-varying topology with energy consumption is considered in this paper. To address the problem of dimension discrepancy among agents, a virtual state is introduced. Based on them, we propose novel PETICs with/without actuation delays to achieve the mean-square exponential consensus of the systems. According to these control schemes, we have achieved the consensus operation of UAVs and UGVs in simulations. In the future work, we will study a class of heterogeneous stochastic MASs with sampled-data outputs and design ETCs based on observers.

\end{document}